\documentclass[12pt, reqno]{amsart}
\usepackage{amsmath}
\usepackage{amssymb}
\usepackage{epsfig}
\usepackage{graphicx}
\usepackage{color}
\definecolor{shadecolor}{gray}{0.875}
\usepackage{amscd}
\usepackage{comment}
\usepackage{adjustbox}
\usepackage{multirow}
\usepackage{tikz-cd}

\input xy
\xyoption{all}

\numberwithin{equation}{section}

\usetikzlibrary{positioning}
\usepackage[colorlinks=false,urlbordercolor=white]{hyperref}
  
\tikzset{sgplattice/.style={inner sep=1pt,norm/.style={red!50!blue},char/.style={blue!50!black},
  lin/.style={black!50}},cnj/.style={black!50,yshift=-2.5pt,left=-1pt of #1,scale=0.5,fill=white}}

\DeclareFontFamily{U}{mathb}{\hyphenchar\font45}
\DeclareFontShape{U}{mathb}{m}{n}{
      <5> <6> <7> <8> <9> <10> gen * mathb
      <10.95> mathb10 <12> <14.4> <17.28> <20.74> <24.88> mathb12
      }{}
\DeclareSymbolFont{mathb}{U}{mathb}{m}{n}
\DeclareMathSymbol{\righttoleftarrow}{3}{mathb}{"FD}

\theoremstyle{plain}
\newtheorem{proposition}{Proposition}[section]
\newtheorem{theorem}[proposition]{Theorem}

\newtheorem{corollary}[proposition]{Corollary}

\newtheorem{lemma}[proposition]{Lemma}

\theoremstyle{definition}
\newtheorem{definition}[proposition]{Definition}

\newtheorem{remark}[proposition]{Remark}

\newtheorem{example}[proposition]{Example}

\newcommand{\eqto}{\stackrel{\lower1.5pt\hbox{$\scriptstyle\sim\,$}}\to}
\newcommand{\eqdashto}{\stackrel{\lower1.5pt\hbox{$\scriptstyle\sim\,$}}\dashrightarrow}

\def\lra{\longrightarrow}

\def\cF{{\mathcal F}}

\def\cX{{\mathcal X}}

\def\fA{{\mathfrak A}}

\def\fS{{\mathfrak S}}

\def\cF{{\overline{\mathcal M}}}

\def\fS{{\mathfrak S}}

\def\bP{{\mathbb P}}

\def\bZ{{\mathbb Z}}
\def\bR{{\mathbb R}}

\def\bC{{\mathbb C}}

\def\rH{{\mathrm H}}

\def\Br{\mathrm{Br}}

\def\CH{\mathrm{CH}}

\def\Pic{\mathrm{Pic}}

\def\Gal{\mathrm{Gal}}
\def\Aut{\mathrm{Aut}}

\def\PGL{\mathsf{PGL}}

\def\lim{\mathrm{lim}}
\def\NS{\mathrm{NS}}

\def\IJ{\mathrm{IJ}}

\def\Cr{\mathrm{Cr}}
\def\cF{\mathcal F}

\author{Tudor Ciurca}
\address{Department of Mathematical Sciences \\
University of Bath \\
Claverton Down \\
Bath \\
BA2 7AY \\
UK}
\email{tc703@bath.ac.uk}

\author{Sho Tanimoto}
\address{Graduate School of Mathematics\\ Nagoya University
\\ Nagoya, 464-8602, Japan}
\email{sho.tanimoto@math.nagoya-u.ac.jp}

\author{Yuri Tschinkel}
\address{Courant Institute\\
                New York University \\
                New York, NY 10012 \\
                USA }
\address{Simons Foundation\\
                 160 Fifth Av.\\ 
                 New York, NY 10010}                
\email{tschinkel@cims.nyu.edu}

\title[Intermediate Jacobians]{Intermediate Jacobians and linearizability}

\begin{document}


\begin{abstract}
We develop an equivariant version of 
the formalism of intermediate Jacobian torsor obstructions, and apply it to conic bundles over rational surfaces, quadric surface bundles over $\bP^1$, and Fano threefolds. 
\end{abstract}

\maketitle

\section{Introduction}
\label{sect:intro}

Let $X$ be a rational variety, over $\bC$, equipped with a generically free, regular action of a finite group $G$. 
A fundamental problem in  higher-dimensional birational geometry is to 
identify {\em linearizable} and {\em projectively linearizable} actions, i.e., actions that are birational to $G$-actions on $\bP(V)$, where  $V$ is a linear representation of $G$, respectively, a linear representation of a central extension of $G$,  see, e.g., \cite[Section 2]{HT-odd} for a discussion of these terms. 
Even the classification of birational types of involutions in dimension 3 is an open problem, see \cite{pro-inv}. 

We develop an equivariant version of 
the formalism of intermediate Jacobian torsor obstructions from \cite{HT-cycle} and \cite{BW19}, and apply it to conic bundles over $\bP^2$, quadric surface bundles over $\bP^1$, and Fano threefolds. 
We pursue the analogy with arithmetic considerations in \cite{HT-quad}, 
\cite{HT-cycle}, \cite{BW19}, \cite{HT-18}, \cite{KuzP}, \cite{FJSV22} and \cite{JJ23}, focused on rationality of geometrically rational threefolds over nonclosed fields.

There are certain similarities between birational geometry over nonclosed fields
and equivariant birational geometry; rationality should be viewed as analogous to linearizability, and birationality to a Brauer-Severi variety as analogous to projective linearizability.  
In the study of rationality, the absolute Galois group of the ground field acts on geometric invariants, such as the Picard group, and all geometric constructions have to take into account Galois symmetries. In the study of (projective) linearizability, the action of the automorphism group limits the choices for birational transformations. 

On the other hand, there are also substantial differences:
\begin{itemize}
    \item Existence of fixed points is not an equivariant birational invariant for actions of nonabelian groups (as can be seen by considering the blowup of the standard 2-dimensional representation of $\fS_3$ at the origin), and  
it is not required for (projective) linearizability (consider the action of $\mathfrak S_3$ on $\bP^1$). 
    \item The action of cyclic groups on projective space is always linearizable (this is an application of the Lefschetz trace formula), while a Galois twist of $\bP^1$ by a cyclic group is not necessarily rational, as is the case for  $\fS_2\simeq\Gal(\bC/\bR)$ acting on a conic over $\bR$ without $\bR$-points.
    \item Rationality of smooth quadric hypersurfaces $X_2\subset \bP^n$, $n\ge 2$, is 
    completely settled, and is equivalent to the existence of rational points. Existence of a $G$-fixed point on a quadric implies linearizability via projection from this point. However, there are also linearizability constructions for quadrics without $G$-fixed points; linearizability is still an open problem for $n\ge 4$, see \cite[Section 9]{TYZ23}.
    \item Some varieties are rational over any field but fail to be linearizable for automorphisms, e.g., a quintic del Pezzo surface or a quintic del Pezzo threefold, which are birationally rigid, and not linearizable, for $G=\fA_5$, see \cite[Section 8.2]{HT23} and \cite{CS}.
    \end{itemize}

\subsection*{Main results and constructions}
In Sections~\ref{sect:jacob} and \ref{sect:equiv-tors} we develop equivariant versions of the theory of intermediate Jacobians and related torsor obstructions. 
Regular actions of finite groups on $X$ yield,
by functoriality, actions on a smooth group scheme $\mathrm{\mathbf{CH}}^2_{X/\bC}$ whose $\bC$-valued points are naturally identified with the Chow group $\CH^2(X)$ of codimension-two cycles. This gives rise to obstructions to projective linearizability analogous to those in \cite[Theorem 6.3]{HT-cycle} and 
\cite[Theorem C]{BW19}:

\begin{theorem}[Theorem~\ref{theorem:IJTobstructions}]
\label{theorem:mainI}
Let $X$ be a smooth projective rational threefold over $\bC$ with a regular, 
projectively linearizable, action of a finite group $G$. 
Then there exists a smooth projective (possibly disconnected) curve $C$ with a regular $G$-action such that for any $G$-invariant connected component $M$ of $\mathrm{\mathbf{CH}}^2_{X/\bC}$  there exist a $G$-invariant connected component $N \subset \Pic(C)$ and an equivariant isomorphism of $G$-varieties
\[
M \cong N.
\]

\end{theorem}

In such situations, we say that the equivariant 
intermediate Jacobian torsor (IJT)
obstruction vanishes. The nonvanishing of this obstruction allows to prove failure of (projective) linearizability 
in many new cases, e.g., for conic bundles 
(see Section~\ref{sect:conic},
Examples \ref{example:involutions_conicbundle}
and \ref{example:C_2timesC_2}), and quadric surface bundles over $\bP^1$ (see Section~\ref{sect:quad}, 
Examples~\ref{example:fixedcurves} and \ref{example:fixedquadric}).  

Furthermore, this gives {\em criteria} for (projective) linearizability, e.g., 
for actions of cyclic groups on conic bundles
$X\to \bP^2$ with quartic discriminant, 
in Theorem~\ref{theorem:rationalitycriterion}, or 
quadric surface bundles $\cX\to \bP^1$ in Theorem~\ref{theorem:quadric_noneffective}.
Among rational Fano threefolds, 
projective linearizability of actions on smooth complete intersections of
two quadrics $X_{2,2}\subset \bP^5$ was settled in \cite{HT-odd}, via a reduction to rationality
considerations over nonclosed fields; 
here we give an alternative proof. We establish a similar criterion for actions on $X_{16}$: they are projectively linearizable if and only if there exists a $G$-invariant rational cubic curve on the variety.  
We  elucidate issues arising for Fano threefolds of type $X_{18}$, see Theorem~\ref{thm:fanolinear}
and the discussion in Section~\ref{sect:fano}. 

Moreover, the intermediate Jacobian torsor formalism yields new general results concerning birationality of (projectively) nonlinearizable actions:

\begin{theorem}[Corollary~\ref{thm:distinguish}]
\label{theorem:main_distinguish}
Let $X_i$, for $i=1,2$, be smooth projective rationally connected threefolds defined over $\bC$ with a regular action of a finite group $G$. Assume that their intermediate Jacobians are the Jacobians 
of smooth projective curves $C_i$,
with transitive action of $G$ on their connected components. Assume that
\begin{itemize}
\item the connected components of $C_i$ have genus $\ge 2$, for $i=1,2$, 
    \item the equivariant IJT obstruction does not vanish for $X_1$, 
    \item  $C_1$ is not $G$-equivariantly isomorphic to $C_2$.
\end{itemize}
Then $X_1$ is not $G$-equivariantly birational to $X_2$.
\end{theorem}

In particular, this applies to {\em involutions} acting on rational conic bundles over rational surfaces or quadric surface bundles over $\bP^1$. Recall that the classification of conjugacy classes of involutions in the Cremona group $\Cr_2(\bC)$, the group of birational automorphisms of $\bP^2$, is based on the study of {\em non-uniruled divisors} in the fixed locus of the involution, see, e.g., \cite{blanc2}. 
In higher dimensions, birational involutions are much more difficult to control: there is too much flexibility.   
The paper \cite{pro-inv} gives a rough classifications of involutions in $\Cr_3(\bC)$, in presence of {\em non-uniruled} divisorial components in the fixed locus. Even then, some of the cases listed in  \cite[Theorem 1.2]{pro-inv} may overlap. In Examples~\ref{example:involutions_conicbundle}, \ref{example:C_2timesC_2}, ~\ref{example:fixedcurves} and \ref{example:fixedquadric}, no such divisors exist; to the best of our knowledge, these
are the first instances  when one is able to distinguish conjugacy classes of such involutions.

An arithmetic version of Theorem~\ref{theorem:main_distinguish} gives many nonrational, geometrically rational, threefolds which are not birational to each other over the ground field; e.g., intersections of two quadrics $X_{2,2}$ considered in \cite{HT-cycle} and \cite{BW19}, or conic bundles in \cite{FJSV22}.  

\

\noindent
{\bf Acknowledgments:}
We are grateful to Chr. B\"ohning, I. Cheltsov, B. Hassett, B. Viray, and O. Wittenberg for comments and suggestions. After the completion of this paper, O. Wittenberg informed us about related results obtained in discussions of F. Scavia with himself, and recorded in unpublished notes, in July 2021. F. Scavia kindly provided these notes \cite{Scavia}, they contain Theorem~\ref{theorem:mainI} and its application to complete intersections of two quadrics. The second author would like to thank Mitsuru Sugata for helpful discussions. Finally, we would like to thank the referee for detailed comments which significantly improved the exposition of the paper.

The first author was partially supported by the Santander Mobility Award.
The second author was partially supported by JST FOREST program Grant number JPMJFR212Z, by JSPS KAKENHI Grand-in-Aid (B) 23K25764, and by JSPS Bilateral Joint Research Projects Grant number JPJSBP120219935.
The third author was partially supported by NSF grant 2301983.

\section{Equivariant intermediate Jacobians}
\label{sect:jacob}

Throughout, we work over the complex numbers $\bC$. 
Let $X$ be a smooth projective rationally connected threefold, $\NS(X)$ its Neron-Severi group, 
and 
\[
\mathrm{IJ}(X) := \rH^3(X, \mathbb C)/(\rH^1(X, \Omega_X^2) \oplus \rH^3(X, \mathbb Z))
\]
its intermediate Jacobian.
This complex torus carries a principal polarization $\theta_X$ induced by the cup product
\[
\wedge^2 \, \rH^3(X, \mathbb Z) \to \rH^6(X, \mathbb Z) \cong \mathbb Z,
\]
so that $(\mathrm{IJ}(X), \theta_X)$ is a principally polarized abelian variety, and we can consider $\mathrm{IJ}(X)$ as a smooth projective variety over $\mathbb C$.

Let $\CH^2(X)$ be the group of codimension-two cycles on $X$, over $\bC$, up to rational equivalence, and $\CH^2(X)_{\rm alg}\subset \CH^2(X)$ the subgroup of cycles algebraically equivalent to $0$. 
Put 
$$
\NS^2(X):=\CH^2(X)/\CH^2(X)_{\rm alg},
$$
it is a finitely generated abelian group, see
\cite[Theorem 1]{BS}.

Let $A$ be an abelian variety over $\mathbb C$ and 
$$
\phi : \CH^2(X)_{\rm alg} \to A(\mathbb C)
$$ 
a group homomorphism. It is called {\em regular} if for any smooth connected variety $T$ over $\bC$, $t_0\in T(\mathbb C)$, and codimension-two cycle $Z \in Z^2(T\times X)$, the map
$$
\begin{array}{rcl}
T(\mathbb C) & \to &  A(\mathbb C), \\
t            & \mapsto &  \phi(Z_t - Z_{t_0}),
\end{array}
$$
is induced by a morphism $T \to A$, defined over $\mathbb C$.
Such a morphism is unique if it exists.
An example of a regular homomorphism is the Abel-Jacobi map
\[
\mathrm{AJ}: \CH^2(X)_{\rm alg} \to \mathrm{IJ}(X)(\mathbb C).
\]
This is bijective by \cite[Theorem I(i)]{BS}, and {\em universal}, i.e., initial object in the category of regular homomorphisms from $\CH^2(X)_{\rm alg}$, see  \cite{Murre}, \cite{Kahn}.  
This identification endows $\CH^2(X)_{\rm alg}$ with the structure of a principally polarized abelian variety, and we denote this scheme by $(\mathrm{\mathbf{CH}}^2_{X/\mathbb C})^0$.

Let $\CH^2(X)^\gamma$ be the preimage of $\gamma \in \NS^2(X)$ under the map 
$$
\CH^2(X) \to \NS^2(X).
$$
The preimage $\CH^2(X)^\gamma$ is in bijection with $\CH^2(X)_{\rm alg}$ (via translation by an element of $\CH^2(X)^\gamma$), which is unique up to a unique translation by elements in $\CH^2(X)_{\rm alg}$. This defines a scheme structure on $\CH^2(X)^\gamma$, and a group scheme structure on $\CH^2(X)$, denoted by  $(\mathrm{\mathbf{CH}}^2_{X/\mathbb C})^\gamma$ and $\mathrm{\mathbf{CH}}^2_{X/\mathbb C}$, respectively. When $X$ is rational, $\mathrm{\mathbf{CH}}^2_{X/\mathbb C}$
is the  group scheme 
that represents the functor 
 $
 \CH^2_{X/\mathbb C, \mathrm{fppf}}
 $ 
 of Chow groups of codimension-two cycles on $X$, 
 constructed in \cite[Theorem 3.1]{BW19}.

A regular action of a finite group $G$ on $X$ induces a $G$-action
on algebraic cycles, and on 
$\CH^2(X)$. It also yields an action on cohomology and a
regular action on the intermediate Jacobian $\mathrm{IJ}(X)$. 
Compatibility with the cup product implies that the polarization $\theta_X\in \NS(\IJ(X))$ is $G$-invariant.

\begin{lemma}
Let $X$ be a smooth projective rationally connected threefold with a regular $G$-action. Then the Abel-Jacobi map
\[
\mathrm{AJ} : \mathrm{CH}^2(X)_{\mathrm{alg}} \to \mathrm{IJ}(X)(\mathbb C)
\]
is $G$-equivariant.
\end{lemma}

\begin{proof}
We recall the construction: let $\gamma$ be a codimension-two cycle which is algebraically equivalent to $0$. 
Algebraic equivalence and homological equivalence coincide, 
provided the Chow group of zero cycles on $X$ is supported on a surface
\cite[Theorem I(ii)]{BS}, which is the case, by the assumption.
Thus $\gamma$ is homologically equivalent to $0$. Let $\alpha$ be a topological $3$-cycle with boundary $\partial \alpha=\gamma$.
    Then the linear functional
    \[
    \int_\alpha  \cdots \in (\rH^1(X, \Omega_X^2)^\vee/\rH_3(X, \mathbb Z)) \cong \mathrm{IJ(X)},
    \]
    is well-defined and is the image of $\alpha$ under the Abel-Jacobi map.
    For any $g \in G$, we have $g\gamma = g\partial \alpha = \partial g\alpha$ on the level of topological cycles. 
\end{proof}

In particular, the $G$-action on $\CH^2(X)_{\rm alg}$ is induced by a regular $G$-action on $(\mathrm{\mathbf{CH}}^2_{X/\mathbb C})^0$. 
\begin{lemma}
    The $G$-action on 
    $
    \CH^2(X)
    $    is induced by a regular $G$-action on $\mathrm{\mathbf{CH}}^2_{X/\bC}$.
\end{lemma}

\begin{proof}
    We  need to verify that for $g \in G$ and $\gamma \in \NS^2(X)$, the action
    \[
    g : (\mathrm{\mathbf{CH}}^2_{X/\bC})^\gamma(\bC) \to (\mathrm{\mathbf{CH}}^2_{X/\bC})^{g\gamma}(\bC), \quad  [Z] \mapsto [gZ],
    \]
    is algebraic. Let $Z_0$ be an effective $1$-cycle representing $\gamma$. The above action is induced by a composition of the morphisms:
    \[
    (\mathrm{\mathbf{CH}}^2_{X/\bC})^\gamma \to (\mathrm{\mathbf{CH}}^2_{X/\bC})^0 \to (\mathrm{\mathbf{CH}}^2_{X/\bC})^0 \to (\mathrm{\mathbf{CH}}^2_{X/\bC})^{g\gamma},
    \]
    mapping
    \[
    [Z] \mapsto [Z]-[Z_0] \mapsto [gZ] - [g Z_0] \mapsto  [gZ],
    \]
    where the middle map is the $G$-action on $(\mathrm{\mathbf{CH}}^2_{X/\bC})^0$ and the other two are translations. 
\end{proof}


A $G$-action on $X$ induces a natural $G$-action on $\mathrm{Chow}^2(X)$, the Chow variety of codimension-two cycles, 
via $Z\mapsto gZ$, for $g \in G$.

\begin{theorem}
\label{thm:tors}
    Let $X$ be a smooth projective rationally connected threefold with a regular action by a finite group $G$.
    The Abel-Jacobi map  
    \[
    \mathrm{AJ} : \mathrm{Chow}^2(X) \rightarrow \mathrm{\mathbf{CH}}^2_{X/\bC}
    \]
    is a $G$-equivariant morphism.
\end{theorem}

\begin{proof}
As in, e.g., 
\cite[Proposition 2.5]{HT-cycle},
one verifies that the Abel-Jacobi map is a well-defined morphism. 
The $G$-equivariance of 
\[
\mathrm{AJ} : \mathrm{Chow}^2(X)(\bC) \to \CH^2(X), \quad Z \mapsto [Z],
\]
is clear from the construction.
\end{proof}

\section{Equivariant intermediate Jacobian torsors and obstructions}
\label{sect:equiv-tors}

\subsection*{Abelian varieties with group actions}
We continue to work over $\bC$. A $G$-abelian variety is an abelian variety with an action of a finite group $G$, preserving the origin.  
Let $A$ be a $G$-abelian variety. 
A principal polarization $\theta_A \in \mathrm{NS}(A)$ is called $G$-equivariant if it is $G$-invariant with respect to the $G$-action on $\mathrm{NS}(A)$; in this case, $(A, \theta_A)$ is called a $G$-equivariant principally polarized abelian variety. 
A $G$-equivariant homomorphism 
$$
\iota: (A, \theta_A)\to (B, \theta_B)
$$
of $G$-equivariant principally polarized abelian varieties is a homomorphism  $\iota : A \to B$  of abelian varieties  such that $\iota^*\theta_B = \theta_A$.
We record:

\begin{lemma}
\label{lemm:decomp}
    Let 
    $$
    \iota: (A, \theta_A)\to (B, \theta_B)
    $$
    be 
a $G$-equivariant homomorphism of $G$-equivariant principally polarized abelian varieties, 
which is a closed embedding. 
Then there exists a $G$-equivariant principally polarized abelian variety $(A', \theta_{A'})$ such that 
$$
(B, \theta_B)\simeq (A\times A', p_1^*\theta_A + p_2^*\theta_{A'}),
$$ 
as $G$-equivariant principally polarized abelian varieties, where $p_i$'s are natural projections.
\end{lemma}

\begin{proof}
    Consider the induced primitive, $G$-equivariant, embedding 
    \[
    \iota_* : \rH_1(A, \mathbb Z) \to \rH_1(B, \mathbb Z).
    \]
    The principal polarization $\theta_B \in \rH^2(B, \mathbb Z) = \wedge^2 \rH^1(B, \mathbb Z)$ defines an alternating form on $\rH_1(B, \mathbb Z)$, and its restriction to $\rH_1(A, \mathbb Z)$ is the alternating form induced by $\theta_A$. This alternating form is compatible with the $G$-action on these lattices. Let $\Lambda$ be the orthogonal complement of $\rH_1(A, \mathbb Z)$ in $\rH_1(B, \mathbb Z)$. Then 
    \[
    \rH_1(A, \mathbb Z) \oplus \Lambda = \rH_1(B, \mathbb Z),
    \]
    since $\theta_B$ and $\theta_A$ are principal polarizations. Since the  alternating form is compatible with $G$-action, $\Lambda$ is $G$-invariant and $\theta_B$ defines a non-degenerate alternating form on $\Lambda$ which is selfdual. 

    The groups $\rH_1(A, \mathbb Z)$ and $\rH_1(B, \mathbb Z)$ admit integral Hodge structures of weight $-1$, which are compatible with 
  $\iota_*$ and preserved by the $G$-action. Since $\theta_B$ is a Hodge class, this induces a principally polarized integral Hodge structure of weight $-1$ on $\Lambda$ which is compatible with $G$-action on it, and  $\Lambda$ can be identified with $\rH_1(A', \mathbb Z)$, where $A'$ is a $G$-equivariant principally polarized abelian variety. 
\end{proof}

\begin{corollary}
\label{coro:decomposition}
    A $G$-equivariant principally polarized abelian variety admits a  decomposition as a product of indecomposable $G$-equivariant principally polarized abelian varieties, which is unique, up to permutation of the factors.
\end{corollary}

An indecomposable $G$-equivariant principally polarized abelian variety is also called {\em irreducible}.

\subsection*{Intermediate Jacobian torsor obstructions}
The following theorem gives an obstruction to equivariant projective linearizability; it is analogous to \cite[Theorem 6.3]{HT-cycle} and 
\cite[Theorem C]{BW19}:

\begin{theorem}
\label{theorem:IJTobstructions}
    Let $X$ be a smooth projective rational threefold with a regular, projectively linearizable, $G$-action.
    Then there exists a smooth projective (possibly disconnected) curve $C$ with a regular $G$-action such that for any $\gamma \in \NS^2(X)^G$, one has a $G$-equivariant isomorphism of $G$-varieties
    $$
    (\mathrm{\mathbf{CH}}^2_{X/\mathbb C})^\gamma\simeq \mathrm{Pic}^m(C),
    $$ 
    for some $G$-invariant class $m$.
    Moreover, we have 
    $$
    (\mathrm{\mathbf{CH}}^2_{X/\mathbb C})^0\simeq \mathrm{Pic}^0(C),
    $$ 
    as $G$-equivariant principally polarized abelian varieties.
\end{theorem}

\begin{proof}
    We have a $G$-equivariant birational map
    \[
    \phi : X \dashrightarrow \mathbb P^3,
    \]
    where $\mathbb P^3$ is equipped with a regular action of $G$. Resolving the indeterminacy of $\phi$ equivariantly, we have equivariant birational morphisms 
    $$
    X \stackrel{g}{\leftarrow} W \stackrel{f}{\rightarrow} \mathbb P^3,
    $$ 
    where $W$ is a smooth projective threefold with a regular $G$-action. By functorial weak factorization (see \cite[Theorem 1.3.3 (1)]{AT19}), $f$ and $g$ are  compositions of equivariant blowups and downs with $G$-irreducible smooth centers. Here $G$-irreducible means that the action on connected components of the smooth center is transitive. The assertion is reduced to a single $G$-equivariant blowup
    \[
    \psi : X_1 \rightarrow X_2,
    \]
    where one of the $X_i$'s satisfies the statement of the assertion, and we need to prove the statement for the other $X_j$.

    Suppose that $X_2$ satisfies the statement, i.e., there exists a smooth projective $G$-curve $C$ such that for any $\gamma \in \NS^2(X_2)^G$, $$(\mathrm{\mathbf{CH}}^2_{X_2/\bC})^\gamma\cong \mathrm{Pic}^m(C),
    $$ 
    for some $m$, as $G$-equivariant varieties, and when $\gamma = 0$, 
    $$
    (\mathrm{\mathbf{CH}}^2_{X_2/\mathbb C})^0\cong \mathrm{Pic}^0(C),
    $$ 
    as $G$-equivariant principally polarized abelian varieties. Suppose that the center of $\psi$ is a $G$-invariant and $G$-irreducible curve $\Gamma$. The blowup formula
    \begin{equation} \label{eqn:decomp}
    \mathrm{CH}^2(X_1) = \mathrm{CH}^2(X_2) \oplus \mathrm{CH}^1(\Gamma),
    \end{equation}
    is $G$-equivariant, and $\psi$ induces a $G$-equivariant isomorphism
    \[
    \rH^3(X_1, \mathbb Z) \cong \rH^3(X_2, \mathbb Z) \oplus \rH^1(\Gamma, \mathbb Z(-1)),
    \]
    compatible with the cup product as well as the Hodge structures. This induces an isomorphism 
    \[
    \mathrm{IJ}(X_1) \cong \mathrm{IJ}(X_2) \times \mathrm{Pic}^0(\Gamma),
    \]
    as principally polarized abelian varieties, and an  isomorphism
    \begin{equation} \label{eqn:iso-ind}
    (\mathrm{\mathbf{CH}}^2_{X_1/\bC})^0\cong \mathrm{Pic}^0(C \sqcup \Gamma),
    \end{equation}
    as $G$-equivariant principally polarized abelian varieties.

   For any $\gamma \in \NS^2(X_1)^G$, 
   \eqref{eqn:decomp}
   induces a $G$-equivariant isomorphism of varieties
   $$
   (\mathrm{\mathbf{CH}}^2_{X_1/\bC})^\gamma \simeq \mathrm{Pic}^m(C \sqcup \Gamma),
   $$ 
   for some $m$, which is a translation of (\ref{eqn:iso-ind}), as claimed.

    When the center of $\psi$ is $0$-dimensional, there is a similar blowup formula for $\mathrm{CH}^2$, after replacing $\mathrm{CH}^1(\Gamma)$ by a permutation module of $\mathbb Z$. In particular, each connected component of $\mathrm{\mathbf{CH}}^2_{X_2/\bC}$ is identified with a connected component of $\mathrm{\mathbf{CH}}^2_{X_1/\bC}$.
    
    In the other direction, when $X_1$ satisfies the assertion, the proof is reversed. Assume that $\psi$ is a blowup of a $G$-invariant and $G$-irreducible curve $\Gamma$. Then we have a $G$-equivariant blowup formula:
    \[
    \mathrm{CH}^2(X_1) = \mathrm{CH}^2(X_2) \oplus \mathrm{CH}^1(\Gamma).
    \]
    As before, $\mathrm{Pic}^0(C)\simeq (\mathrm{\mathbf{CH}}^2_{X_2/\bC})^0 \times \mathrm{Pic}^0(\Gamma)$ as $G$-equivariant principally polarized abelian varieties.
    When every connected component of $\Gamma$ has genus $0$, the assertion is clear.
    Thus, we may assume that every connected component of $\Gamma$ has genus $g \geq 1$.
    Using Corollary~\ref{coro:decomposition} and the  Torelli theorem for curves, we can write $C = C'\sqcup \Gamma$, where $C'$ is a smooth projective $G$-curve, and we have
    \[
    (\mathrm{\mathbf{CH}}^2_{X_2/\bC})^0\cong \mathrm{Pic}^0(C').
    \]
    Then, for any $\gamma \in \NS^2(X_2)^G$, $$(\mathrm{\mathbf{CH}}^2_{X_2/\bC})^\gamma \oplus \mathrm{Pic}^0(\Gamma)\cong \mathrm{Pic}^m(C),
    $$
    for some $m$, as $G$-varieties. 
    Since an isomorphism
    $$
    (\mathrm{\mathbf{CH}}^2_{X_2/\bC})^\gamma \oplus \mathrm{Pic}^0(\Gamma)\cong \mathrm{Pic}^m(C) = \mathrm{Pic}^{m'}(C') \oplus \mathrm{Pic}^{0}(\Gamma)
    $$ 
    is a translation of an isomorphism $$(\mathrm{\mathbf{CH}}^2_{X_2/\bC})^0 \oplus \mathrm{Pic}^0(\Gamma)\cong \mathrm{Pic}^0(C) = \mathrm{Pic}^{0}(C') \oplus \mathrm{Pic}^{0}(\Gamma),
    $$ 
    and the translation preserves the direct sum structure, we see that 
    $$
    (\mathrm{\mathbf{CH}}^2_{X_2/\bC})^\gamma\cong \mathrm{Pic}^{m'}(C'),
    $$ 
    as $G$-varieties. 
When the center of $\psi$ is $0$-dimensional, the proof is similar, and we omit it. 
\end{proof}

The intermediate Jacobian torsor formalism yields new obstructions to equivariant birationality of nonlinearizable actions on rationally connected threefolds, and is a new tool to distinguish finite subgroups of the Cremona group $\Cr_3(\bC)$, up to conjugation in this group. 
Before stating the main theorem, let us introduce one terminology:
\begin{definition}
    Let $X$ be a smooth projective rationally connected threefold with a regular $G$-action. We say that {\em the equivariant {\rm IJT}-obstruction} vanishes for $X$ if there exists a smooth projective (possibly disconnected) curve $C$ with a regular $G$-action such that for any $\gamma \in \NS^2(X)^G$, one has a $G$-equivariant isomorphism of $G$-varieties
    $$
    (\mathrm{\mathbf{CH}}^2_{X/\mathbb C})^\gamma\simeq \mathrm{Pic}^m(C),
    $$ 
    for some $G$-invariant class $m$.
    Moreover, we have 
    $$
    (\mathrm{\mathbf{CH}}^2_{X/\mathbb C})^0\simeq \mathrm{Pic}^0(C),
    $$ 
    as $G$-equivariant principally polarized abelian varieties.
\end{definition}

As an application of Theorem~\ref{theorem:IJTobstructions}, we obtain:

\begin{theorem}
\label{theorem:distinguish}   
Let $X_i$, with $i=1,2$, be smooth projective rationally connected $G$-varieties of dimension 3 such that 
their intermediate Jacobians are irreducible $G$-equivariant principally polarized abelian varieties.
Assume that 
\begin{itemize}
    \item the equivariant {\rm IJT}-obstruction does not vanish for $X_1$,
    \item $\IJ(X_1)$ is not isomorphic to $\IJ(X_2)$ as $G$-equivariant principally polarized abelian varieties. 
\end{itemize}    
Then $X_1$ is not $G$-equivariantly birational to $X_2$.
\end{theorem}

\begin{proof}
    Suppose that $X_1$ is $G$-equivariantly birational to $X_2$. By assumption, for any smooth projective $G$-curve $C$
 there is a $G$-invariant class $\gamma$ of codimension-two cycles on $X_1$ such that $$
 (\mathrm{\mathbf{CH}}^2(X_1))^\gamma\not \cong \Pic^e(C),
 $$ 
  for any $G$-invariant $e$. As in Theorem~\ref{theorem:IJTobstructions}, this shows that $\mathrm{IJ}(X_1)$ must arise in the decomposition of $\mathrm{IJ}(X_2)$.
    Since they are irreducible, as $G$-principally polarized abelian varieties, this implies that $\IJ(X_1) \cong \IJ(X_2)$. This contradicts the assumption. 
\end{proof}

\begin{corollary}
\label{thm:distinguish}   
Let $X_i$, with $i=1,2$, be smooth projective rationally connected $G$-varieties of dimension 3 such that 
their equivariant intermediate Jacobians are the Jacobians of $G$-irreducible smooth projective curves $C_i$ whose connected components have genus $g_i \geq 2$. 
Assume that 
\begin{itemize}
    \item the equivariant {\rm IJT}-obstruction does not vanish for $X_1$,
    \item $C_1$ is not $G$-equivariantly isomorphic to $C_2$. 
\end{itemize}    
Then $X_1$ is not $G$-equivariantly birational to $X_2$.
\end{corollary}
\begin{proof}
By Torelli theorem (see, e.g., \cite[Theorem 12.1]{Milne}), $\Pic^0(C_1)$ and $\Pic^0(C_2)$ are not isomorphic to each other as $G$-equivariant principally polarized abelian varieties.  The assertion follows from Theorem~\ref{theorem:distinguish}.
\end{proof}

\section{Equivariant Prym varieties}
\label{subsec:equivariantPrym}

We consider $G$-equivariant Prym varieties, 
following the constructions in \cite[Section 4]{FJSV22} in the arithmetic context. 

\subsection*{Equivariant Prym varieties}

Fix a
$G$-equivariant \'etale finite morphism 
$\varpi : \widetilde{\Delta} \to \Delta$ 
between connected smooth projective $G$-curves, of degree $2$.
The norm map between the Picard schemes
\[
\mathrm{Nm} : \Pic(\widetilde{\Delta}) \to \Pic(\Delta),
\]
is a $G$-equivariant homomorphism of group schemes.
A $G$-equivariant embedding $r: \Delta \hookrightarrow W$ 
into a smooth projective rational $G$-surface induces a $G$-equivariant morphism of group schemes
\[
r^* : \Pic(W) \to \Pic(\Delta).
\]
The $\Pic(W)$-polarized Prym scheme is defined as
\[
\mathrm{PPrym}^{\Pic(W)}(\widetilde{\Delta}/\Delta) := \Pic(\widetilde{\Delta}) \times_{\Pic(\Delta)} \Pic(W),
\]
it is a group scheme with diagonal $G$-action. The identity component of $\mathrm{PPrym}^{\Pic(W)}(\widetilde{\Delta}/\Delta)$ is the identity component of the kernel of the norm map $\mathrm{Nm}$ which is the Prym variety $\mathrm{Prym}(\widetilde{\Delta}/\Delta)$ from \cite{Beau77}. This is a principally polarized abelian variety, and $G$ acts on it. The polarization associated with a theta divisor on $\Pic^0(\widetilde{\Delta})$ is a $G$-principal polarization. Its restriction to $\mathrm{Prym}(\widetilde{\Delta}/\Delta)$ is twice a principal polarization on $\mathrm{Prym}(\widetilde{\Delta}/\Delta)$. In particular, $\mathrm{Prym}(\widetilde{\Delta}/\Delta)$ admits a $G$-invariant principal polarization.

For $D \in \Pic(W)(\bC)$, denote the fiber of 
$$
\mathrm{PPrym}^{\Pic(W)}(\widetilde{\Delta}/\Delta) \to \Pic(W)
$$ 
above $D$ by $V_D$; it is a $G$-equivariant torsor of $V_0$ when $D$ is $G$-invariant.
Since $\varpi_*\varpi^*D = 2D$, we have 
\[
V_{2D + D'} = \varpi^*D + V_{D'}, \quad \forall D, D' \in \Pic(W)(\bC),
\]
compatibly with $G$-actions. The group scheme $V_0$ consists of two connected components, 
\[
P = \mathrm{Prym}(\widetilde{\Delta}/\Delta), \quad \widetilde{P} = V_0 \setminus P.
\]
Both carry the $G$-action, and $\widetilde{P}$ is an equivariant $2$-torsion torsor of $P$.

For $W = \mathbb P^2$, $\Delta \subset \mathbb P^2$ a $G$-invariant smooth plane quartic, and $H$ the hyperplane class of $\mathbb P^2$ we have:

\begin{proposition}[{\cite[Section 6, Equation (6.1)]{prym}}]
The parity of $h^0$ is constant on each 
of the two connected components of $V_{H}$, and is different on these components.
\end{proposition}

\begin{definition}[{\cite[Definition 4.4]{FJSV22}}]
    Let $P^{(1)}$ be the connected component of $V_H$ on which $h^0$ is even, and let $\widetilde{P}^{(1)}$ be the connected component of $V_H$ on which $h^0$ is odd.  
For each $m\in \bZ_{\ge 0}$ and $e  = 0, 1$, we define
    \[
    P^{(2m + e)} = P^{(e)} + mH, \quad \widetilde{P}^{(2m+e)} = \widetilde{P}^{(e)} + mH,
    \]
    where $P^{(0)} = P$ and $\widetilde{P}^{(0)} = \widetilde{P}$.
    These are $G$-equivariant torsors of $\mathrm{Prym}(\widetilde{\Delta}/\Delta)$.
\end{definition}

\subsection*{Intermediate Jacobians and Prym varieties}

We follow \cite{Beau77} and \cite{FJSV22}. 
Let $X$ be a smooth projective rational threefold over $\bC$ with a regular action of a finite group $G$, and admitting a 
$G$-equivariant standard conic bundle structure $\pi: X\to W$,
where 
\begin{itemize}
    \item  $W$ is a smooth projective rational $G$-surface,
    \item the discriminant curve $\Delta\subset W$ is smooth, 
    \item the (equivariant) \'etale double cover 
    $\varpi : \widetilde{\Delta}\to \Delta$ is irreducible. 
\end{itemize}

\begin{lemma}[{\cite[Proposition 5.3]{FJSV22}}]
    Let $\gamma$ be the algebraic class of a line in a singular fiber of $\pi$. Then $\gamma$ is $G$-invariant and we have a $G$-equivariant exact sequence:
    \[
    0 \to \bZ \gamma \to \mathrm{NS}^2(X) \stackrel{\pi_*}{\lra} \mathrm{NS}(W) \to 0.
    \]
\end{lemma}

\begin{proof}
Since $\pi$ is $G$-equivariant, $\pi_*$ is as well. Moreover, $G$ maps any line in the singular fibers of $\pi$ to a line in a singular fiber. Since $\widetilde{\Delta}$ is irreducible, $\gamma$ is $G$-invariant. The exactness of the sequence follows from \cite[Proposition 5.3]{FJSV22}. 
\end{proof}

The following is an analog of a key result in \cite{FJSV22}, comparing torsors of intermediate Jacobians and Prym varieties:

\begin{theorem}[{\cite[Theorems 5.1 and 5.8]{FJSV22}}]
\label{theorem:ChowvsPrym}
For $\pi : X \to W$ as above, 
 there is a $G$-equivariant surjective morphism of group schemes
    \[
    \mathrm{\mathbf{CH}}^2_{X/\bC} \to \mathrm{PPrym}^{\Pic(W)}(\widetilde{\Delta}/\Delta)
    \]
    that is an isomorphism when restricted to each connected component.
    Moreover, this induces an isomorphism
    \[
    (\mathrm{\mathbf{CH}}^2_{X/\bC})^0 \to \mathrm{Prym}(\widetilde{\Delta}/\Delta)
    \]
    of $G$-equivariant principally polarized abelian varieties.
\end{theorem}

\begin{proof}
     A morphism as stated is established in \cite[Theorem 5.1]{FJSV22}, so we only need to show that it is $G$-equivariant. We recall the constructions  in \cite[Section 5]{FJSV22}:

    We have a $G$-equivariant embedding $\iota : \Delta \hookrightarrow X$ with $G$-invariant image. Let $\epsilon : X' \to X$ be the $G$-equivariant blowup along $\iota(\Delta)$. Let $X_\Delta:=\pi^{-1}(\Delta)$ and $S$ be its proper transform in $X'$, embedded via $j : S \hookrightarrow X'$. Note that $X_\Delta$ is $G$-invariant, and so is $S$. Then $S$ admits a $G$-equivariant 
    $\mathbb P^1$-bundle $p : S \to \widetilde{\Delta}$. We have a homomorphism 
    \[
    p_*j^*\epsilon^* : \mathrm{CH}^2(X) \to \mathrm{CH}^1(\widetilde{\Delta})= \Pic(\widetilde{\Delta}),
    \]
    where $\epsilon^* : \mathrm{CH}^2(X) \to \mathrm{CH}^2(X')$ and $j^* : \mathrm{CH}^2(X') \to \mathrm{CH}^2(S)$ are the refined Gysin homomorphisms from \cite[Section 6.6]{Fulton}. By  \cite[Theorem 6.2]{Fulton}, these are $G$-equivariant, so the above homomorphism is also $G$-equivariant.
 By \cite[Proposition 5.6]{FJSV22}, the homomorphism
    \[
    (p_*j^*\epsilon^*, \pi_*):\mathrm{CH}^2(X) \to \Pic(\widetilde{\Delta}) \times \Pic(W)
    \]
    is $G$-equivariant, with image contained in $\mathrm{PPrym}^{\Pic(W)}(\widetilde{\Delta}/\Delta)$. This is the homomorphism used in \cite[Theorem 5.1]{FJSV22}. The last statement follows from \cite[Lemma 5.5 and Theorem 5.8]{FJSV22}. 
\end{proof}

\section{Conic bundles}
\label{sect:conic}

We turn to the linearizability of actions on smooth rational Fano threefolds $X$ over $\bC$ admitting an equivariant conic bundle structure 
$$
\pi : X \to \mathbb{P}^2,
$$
such that the relative Picard rank of $\pi$ is $1$ and the {\em discriminant curve} $\Delta\subset \bP^2$ is a smooth plane quartic. We develop a version of intermediate Jacobian torsor obstructions in this setting, following \cite{FJSV22}, and establish a criterion for linearizability of actions of cyclic groups.

\subsection*{Classification}
The classification of Fano threefolds
implies that  $X$ is one of the following, see 
\cite{Mori1981}, \cite{MM83}, \cite{MM03}, and \cite{Matsuki}:
\begin{enumerate}
\item a double cover of $\mathbb{P}^1 \times \mathbb{P}^2$ branched along a smooth divisor of type $(2, 2)$,
\item the blowup of a quadric hypersurface in $\mathbb{P}^4$ along a genus $2$, degree $6$ curve,  
\item the blowup of a smooth intersection of two quadrics in $\mathbb{P}^5$ along a plane conic. 
\end{enumerate}
In all cases, $X$ is rational, see, e.g., \cite[Corollary 5.6.1]{Prokhorov2018}, and the automorphism group 
$
\Aut(X)
$ 
is finite. 

A $G$-action on $X$ induces a $G$-action on the Picard group $\Pic(X)$, which has 
to preserve the extremal rays of the nef cone.
In Cases (2)-(3), this implies that the $G$-action necessarily preserves the blowup presentation. Thus, the linearization problem is reduced to the case of $G$-actions on quadrics (considered in \cite{TYZ23}) or 
the case of intersections of two quadrics (considered in \cite{HT23}). Linearizability is still largely open in Case (2), see \cite[Section 9]{TYZ23}. In Case (3), projective linearizability is equivalent to the existence of a $G$-stable line.

\subsection*{The main example}
Thus, we focus on Case (1), the double cover
\begin{equation} \label{eqn:X}
\delta: X\to \bP^1\times \bP^2,
\end{equation}
ramified in a smooth divisor of bi-degree $(2,2)$, with projections

\centerline{
\xymatrix{
X\ar[r]^{\pi_2} \ar[d]_{\pi_1}&  \bP^2 \\
 \bP^1 & 
}
}
\noindent
where $\pi_2$ is
a standard conic bundle, with discriminant a smooth quartic  $\Delta\subset \bP^2$ and discriminant 
double cover $\varpi : \widetilde{\Delta} \to \Delta$.
Note that $\Delta$ being smooth does not follow from smoothness of $X$, and we impose this additional assumption throughout this section.

A regular $G$-action 
on $X$ induces actions of $G$ on $\bP^1,\bP^2$, $\Delta$, and $\widetilde{\Delta}$, so that $\pi_1$, $\pi_2$ and $\varpi : \widetilde{\Delta} \to \Delta$ are $G$-equivariant. 
Let $\mathcal F_1(X/\bP^1)$ be the Fano variety of lines in the fibers of $\pi_1$; it 
carries the $G$-action as well, and we have
an equivariant Stein factorization 
\[
\mathcal F_1(X/\bP^1) \to C \to \mathbb P^1,
\]
where $C$ is 
a smooth projective curve of genus 2 equipped with a regular $G$-action.

\subsection*{Standard linearization construction}

Over non-closed fields, $X$ is $k$-rational if $\widetilde{\Delta}(k)\neq \emptyset$. An analogous statement in the equivariant setting is the following proposition:

\begin{proposition}[{\cite{Prokhorov2018}}]
\label{proposition:prokhorov2018}
Assume that $G$ is cyclic and $\widetilde{\Delta}^G\neq \emptyset$. Then the $G$-action on $X$ is linearizable. 
\end{proposition}

\begin{proof}
    Let $H_1$, $H_2$ be the pullbacks of hyperplane classes from $\mathbb P^1$ and $\mathbb P^2$, respectively.
    Let $\ell\cup \ell'$ be a singular fiber of $\pi_2$ corresponding to the $G$-fixed point on $\tilde{\Delta}$, i.e., both $\ell$ and $\ell'$ are $G$-stable.
    Let $p : \widehat{X} \to X$ the blowup of $X$ along $\ell$, with exceptional divisor $E$.
    Then $\widehat{X}$ is a smooth weak Fano threefold with a regular $G$-action induced by the $G$-action on $X$. 

    Let $L = p^*H_1 + p^*H_2 - E$; we have $L^3 = 2$. The linear system $|L|$ is $4$-dimensional and base point free,  the associated morphism 
    $$
    \Phi_{|L|} : \widehat{X} \to \bP^4
    $$ 
    defines a $G$-equivariant birational morphism $\widehat{X} \to Q\subset \bP^4$
    to a quadric threefold  
    which is a quadric cone of rank $4$; with linear $G$-action on the ambient $\bP^4$. 
    Now $\ell$ and $\ell'$ meet at a $G$-fixed point. The fiber $S$ of $\pi_1 : X \to \mathbb P^1$ containing this fixed point is $G$-stable. Let $\widehat{S}$ be its strict transform on $\widehat{X}$. Then $L^2.\widehat{S} = 1$, i.e., $Q$ contains a $G$-stable plane. Since $G$ is cyclic, $Q$ contains a $G$-fixed smooth point. Indeed, a $G$-stable plane contains at least $5$ fixed points by the fixed point formula. Projection from this point yields an equivariant birational map to $\mathbb P^3$. 
    \end{proof}

\subsection*{Intermediate Jacobians and torsors}

\begin{theorem} \cite{Bruin}, \cite[Theorem 4.5]{FJSV22}
\label{theorem:Prymvariety}
    We have equivariant isomorphisms
    \[
    P^{(0)} = \mathrm{Prym}(\widetilde{\Delta}/\Delta) \cong \Pic^0(C), \quad P^{(1)} \cong \Pic^1(C).
    \]
    Moreover, the first isomorphism is an isomorphism of $G$-equivariant principally polarized abelian varieties.
\end{theorem}

\begin{proof}
   Choose a general $x\in C$ with image $t \in \mathbb P^1$. The fiber $S_t$ of $\pi_1$ is a smooth quadric surface and $x$ corresponds to a ruling of lines on $S_t$. Let $\ell$ be a general line in this ruling. The pushfoward of $\ell$ to $\mathbb P^2$ via $\pi_2$ is a line so that  $\pi_2^*\Delta .\ell = 4$. Since $x$ and $\ell$ are general, $\ell$ meets  $\pi_2^*\Delta$ at smooth points of $\pi_2^*\Delta$. In particular, $\pi_2^*\Delta\cap \ell$ and its pushforward to $\widetilde{\Delta}$ define a degree $4$ effective $0$-cycle $\alpha$ on $\widetilde{\Delta}$. Since $\ell$ is parametrized by $\mathbb P^1$, the linear equivalence class of $\alpha$ only depends on $x$, and we have a natural $G$-equivariant morphism
   \[
   C \to \Pic(\widetilde{\Delta}).
   \]
   Since the pushfoward of $\alpha$ to $\Delta$ is the intersection of $\Delta$ and a line on $\mathbb P^2$, the above morphism actually defines
   \[
   C \to V_{K_\Delta}.
   \]
   Since $\alpha$ is parametrized by $\mathbb P^1$, we have $h^0(\widetilde{\Delta}, \mathcal O(\alpha)) \geq 2$. By Clifford's theorem, $h^0(\widetilde{\Delta}, \mathcal O(\alpha)) = 2$.
   Hence we have a $G$-equivariant morphism
   \[
   C \to P^{(1)},
   \]
   inducing a $G$-equivariant morphism
   \[
   \Pic^1(C) \to P^{(1)},
   \]
which is an isomorphism by \cite{Bruin} and \cite[Theorem 4.5]{FJSV22}. The above construction of this morphism is slightly different from the one in \cite{Bruin} and \cite[Theorem 4.5]{FJSV22}, however they are compatible and they define the same morphism. (See \cite[Proposition 6.3(iii)]{FJSV22} for more details.) 
Finally, the above $G$-equivariant isomorphism induces a $G$-equivariant isomorphism
   \[
   \Pic^0(C) \to P^{(0)}.
   \]
  For the last claim see, e.g., \cite[Section 5, Case 4]{Bruin}. 
\end{proof}

Let $\gamma_1$ be the algebraic class of lines in the fibers of $\pi_1:X\to \mathbb P^1$ and 
$\gamma_2$ the algebraic class of lines in the fibers of 
$\pi_2: X\to \mathbb P^2$.  Then $\NS^2(X)$ is generated by $\gamma_1$ and $\gamma_2$ and both classes are $G$-invariant.

\begin{theorem}[{\cite[Theorem 6.4]{FJSV22}}]
\label{theorem:FJSV226.4}
    In the above setting, we have equivariant isomorphisms:
    \[
    (\mathrm{\mathbf{CH}}^2_{X/\bC})^{m\gamma_1+n\gamma_2} \cong
    \begin{cases}
      P & \text{$m, n$ are even;}\\
      \widetilde{P} & \text{$m$ is even and $n$ is odd;}\\
      P^{(1)} & \text{$m$ is odd and $n$ is even;}\\
      \widetilde{P}^{(1)} & \text{$m, n$ are odd.}
    \end{cases}
    \]
\end{theorem}
\begin{proof}
The analogous isomorphisms in \cite[Theorem 6.4]{FJSV22} are induced by the morphism in Theorem~\ref{theorem:ChowvsPrym}, which is 
equivariant.
\end{proof}

\subsection*{Equivariant IJT obstructions}

The following is the equivariant analog of \cite[Theorem 1.5]{FJSV22}:

\begin{theorem}
\label{theorem:rationalitycriterion}
Let $G$ be a cyclic group and 
$\delta: X\to \mathbb P^1 \times \mathbb P^2$ a smooth $G$-double cover ramified along a smooth divisor of bi-degree $(2, 2)$.
Assume that $\Delta$ is smooth.
Then $X$ is $G$-linearizable if and only if $\widetilde{P}$ or $\widetilde{P}^{(1)}$ is a trivial torsor.
\end{theorem}

\begin{remark}
The triviality of $\Br(k)$ in the assumptions in \cite[Theorem 1.5]{FJSV22}    
translates into the condition $\rH^2(G,k^\times) =0$.  
\end{remark}

We start with auxiliary constructions.   
Consider
\[
\overline{M}_{0,0}(X, \gamma), \quad \gamma \in \NS^2(X), 
\]
the coarse moduli space of stable maps of genus $0$ and class $\gamma$ and put
$$
\cF_{1,1}(X):= \overline{M}_{0,0}(X, \gamma_1+\gamma_2).
$$

\begin{proposition}
\label{proposition:N11}
The moduli space $\cF_{1, 1}(X)$ is smooth and the Abel-Jacobi map 
$$
\mathrm{AJ} : \cF_{1, 1}(X) \to (\mathrm{\mathbf{CH}}^2_{X/\mathbb C})^{\gamma_1 + \gamma_2}
$$ 
is a $\mathbb P^1$-fibration.
\end{proposition}

\begin{proof}
    There are two types of stable maps parameterized by $\cF_{1,1}(X)$:
    \begin{itemize}
        \item[(1)] an embedding $f : \bP^1 \to X$ of class $\gamma_1 + \gamma_2$;
        \item[(2)] the domain consists of two $\bP^1$'s, one is isomorphically mapping to a line $R_1$ in a fiber of $\pi_1$ and the other to a line $R_2$ in a fiber of $\pi_2$.
    \end{itemize}
    These have trivial automorphisms. Thus, for the first claim, we only need to show the smoothness of the moduli stack 
    \[
    \overline{\mathcal M}_{0,0}(X, \gamma_1 + \gamma_2).
    \]

    Consider (1), an embedding $f : R \cong \bP^1 \to X$, and let $\ell:=\pi_2(f(R))$ be its image in $\mathbb P^2$, a line.
    Its preimage $S_\ell:= \pi_2^{-1}(\ell)$ is a normal surface with only canonical singularities. Consider the exact sequences of tangent and normal bundles
    \[
    0 \to T_R \to T_{S_\ell}|_R \to Q \to 0,
    \]
    and
    \[
    0 \to Q \to N_{R/X} \to N_{S_\ell/X}|_R \to 0.
    \]
    Since $N_{R/X}$ is locally free, $Q$ is torsion-free so that it is locally free. Thus, $T_{S_\ell}|_R$ is also torsion-free.
    We have 
    \[
    -K_{S_\ell}.R = (-K_X - S_{\ell}).R = 3-1 = 2
    \]
    This indicates that $c_1(Q).R = 0$ and $Q \cong \mathcal O$. We also have $S_\ell .R = 1$, so that $N_{S_\ell/X}|_R\cong \mathcal O(1)$ and $N_{R/X} \cong \mathcal O \oplus \mathcal O(1)$; we have $\rH^1(R, N_{R/X}) = 0$.
    Thus $[f]$ is a smooth point of $
    \mathcal F_{1, 1}(X)$.

    Consider an $f : R =  R_1 \cup R_2 \to X$, from Case (2).
    As above, 
    $$
    N_{R_1/X}\cong \mathcal O \oplus \mathcal O \quad \text{ and } \quad N_{R_2/X} \cong \mathcal O \oplus \mathcal O(-1),
    $$
    so that $\rH^1(R_i, f^*T_X|_{R_i}) = 0$ for $i = 1, 2$.
    Since 
    $$
    T_X|_{R_1} \cong \mathcal O(2) \oplus \mathcal O \oplus \mathcal O,
    $$
    it follows from \cite[Proposition 2.3]{BLRT23} that
    \[
    \rH^1(R, f^*T_X) = 0.
    \]
    By the discussion of \cite[Section 1.1]{Testa},  $[f]\in \cF_{1, 1}(X)$ is a smooth point.
    Hence $\cF_{1, 1}(X)$ is a smooth projective variety of dimension $3$ so that the Abel-Jacobi map
    \[
     \mathrm{AJ}: \cF_{1, 1}(X) \to (\mathrm{\mathbf{CH}}^2_{X/\mathbb C})^{\gamma_1 + \gamma_2}
    \]
    is well-defined. We proceed to show show that it is
a $\mathbb P^1$-fibration. This shows that $\cF_{1, 1}(X)$ is irreducible.

A stable map $[f : R \to X] \in \cF_{1, 1}(X)$ defines a line $\ell:=\pi_2(f(R))$ in $\mathbb P^2$ and thus 
    a morphism to the dual
    \[
    \cF_{1, 1}(X) \to (\bP^2)^\vee.
    \]
We denote its Stein factorization  by
    \[
    \cF_{1, 1}(X) \to B \to (\bP^2)^\vee.
    \]
    The pullback $S_\ell := \pi_2^{-1}(\ell)$ is a normal projective surface with only canonical singularities, and $f(R)$ is a Weil divisor on $S_\ell$. A fiber of $\cF_{1, 1}(X) \to B$ is the linear system $|f(R)|$ of the divisorial sheaf $\mathcal O_{S_\ell}(f(R))$ which is isomorphic to $\mathbb P^1$. Indeed, a rational curve $[f : R \to S_\ell]$ deforms to cover $S_\ell$ so its moduli has expected dimension $1$. Hence $\cF_{1, 1}(X) \to B$ is a $\bP^1$-fibration, and the Abel-Jacobi map factors as
    \[
    \cF_{1, 1}(X) \to B \to (\mathrm{\mathbf{CH}}^2_{X/\mathbb C})^{\gamma_1 + \gamma_2}.
    \]
    The assertion follows if we show that $\mathrm{AJ}$ 
    is dominant with connected fibers. Recall that the Abel-Jacobi map is induced by
    \[
    \cF_{1, 1}(X) \to \Pic^4(\widetilde{\Delta}),
    \]
    given in the proof of Theorem~\ref{theorem:ChowvsPrym}, which factors as
    \[
    \cF_{1, 1}(X) \to B \to \Pic^4(\widetilde{\Delta}).
    \]
    To show that $B \to \Pic^4(\widetilde{\Delta})$ is injective, 
    fix $[f] \in \cF_{1, 1}(X)$ and let $D$ be its image  in $\Pic^4(\widetilde{\Delta})$. We claim that $h^0(\widetilde{\Delta}, \mathcal O(D)) = 1$. Otherwise, by Clifford's theorem,  $h^0(\widetilde{\Delta}, \mathcal O(D)) = 2$.
    However, by Theorem~\ref{theorem:FJSV226.4}, the image of the Abel-Jacobi map is in $\widetilde{P}^{(1)}$. 
    This contradicts the fact that $h^0(\widetilde{\Delta}, \mathcal O(D))$ is even. Thus we have $h^0(\widetilde{\Delta}, \mathcal O(D)) = 1$. This implies that $B \to \Pic^4(\widetilde{\Delta})$ is injective. 
\end{proof}

The following is analogous to \cite[Proposition 6.5]{FJSV22}:

\begin{proposition}
\label{proposition:11}
    Let $G$ be a cyclic group.
    Assume that $\widetilde{P}^{(1)}$ is a trivial $P$-torsor.
    Then $X$ is $G$-linearizable.
\end{proposition}

\begin{proof}
    Assume that $\widetilde{P}^{(1)}$ has a $G$-fixed point. By Theorem~\ref{theorem:FJSV226.4},  
    \[
(\mathrm{\mathbf{CH}}^2_{X/\mathbb C})^{\gamma_1 + \gamma_2},
    \]
    has a fixed point. Since $G$ is cyclic, Proposition~\ref{proposition:N11} implies that
    \[
    \cF_{1, 1}(X)
    \]
    has at least two fixed points. In particular, $\pi_1 : X \to \mathbb P^1$ admits two $G$-invariant sections.
    Projection from one of  the sections gives a $\bP^2$-bundle over $\mathbb P^1$, with a $G$-fixed section.
    Such a $\bP^2$-bundle is $G$-linearizable.
\end{proof}

Consider the coarse moduli space
\[
\cF_{2, 1}(X) = \overline{M}_{0,0}(X, 2\gamma_1 + \gamma_2).
\]
Dimension count shows that it generically parameterizes a stable map from $\bP^1$ which is an embedding into $X$.

Consider $\mathrm{Sym}^{8}(\widetilde{\Delta}) \to \mathrm{Pic}^8(\widetilde{\Delta})$ and denote the preimage of $\widetilde{P}^{(2)} \subset \mathrm{Pic}^8(\widetilde{\Delta})$ by $\widetilde{S}^{(2)}$. Then 
$$
\widetilde{S}^{(2)} \to \widetilde{P}^{(2)}
$$ 
is a $\mathbb P^3$-bundle over $\widetilde{P}^{(2)}$.
For a general smooth irreducible curve $C$ of type $2\gamma_1 + \gamma_2$, $C \cap \pi_2^{-1}(\Delta)$ induces a degree-$8$ divisor on $\widetilde{\Delta}$ which defines a class on $\widetilde{P}^{(2)}$.
After taking a $G$-equivariant resolution $\widetilde{\cF}_{2, 1}(X)$, we have 
a $G$-equivariant birational morphism
\[
\widetilde{\cF}_{2, 1}(X) \to \widetilde{S}^{(2)}.
\]

\ 

The following is analogous to \cite[Theorem 6.7]{FJSV22}:

\begin{proposition}
\label{proposition:12}
    Assume that $G$ is cyclic and $\widetilde{P}$ is a trivial $P$-torsor.
    Then $X$ is $G$-linearizable.
\end{proposition}

\begin{proof}
If $\widetilde{P}$ is trivial then so is $\widetilde{P}^{(2)}$, and $\widetilde{S}^{(2)}$ has a fixed point. 
Since the existence of fixed points is an equivariant birational invariant, for actions of abelian groups on smooth projective varieties, both $\widetilde{\cF}_{2, 1}(X)$ and $\cF_{2, 1}(X)$ have a fixed point. Thus, there is a $G$-invariant section of $\pi_1 : X \to \bP^1$. Projecting from it gives a $\mathbb P^2$-bundle over $\bP^1$. The generic fiber $X_\eta$ is isomorphic to $\mathbb P^2_\eta$ with a $G$-invariant line. Such bundles are linearizable, for cyclic actions.
\end{proof}

\begin{proof}[Proof of Theorem~\ref{theorem:rationalitycriterion}]
If $X$ is $G$-linearizable, then Theorems~\ref{theorem:IJTobstructions}, ~\ref{theorem:Prymvariety}, and ~\ref{theorem:FJSV226.4} show 
that $\widetilde{P}$ and $\widetilde{P}^{(1)}$ are isomorphic to $P$ or $P^{(1)}$.
Since
$$
[\widetilde{P}^{(1)}] = [P^{(1)}] + [\widetilde{P}] \in \rH^1(G, P),
$$
$\widetilde{P}$ or $\widetilde{P}^{(1)}$ is a trivial torsor.

Conversely, when $\widetilde{P}^{(1)}$, respectively $\widetilde{P}$, is trivial, 
we apply Proposition~\ref{proposition:11},
respectively ~\ref{proposition:12}.

\end{proof}

\subsection*{Examples of (non)linearizable actions}
 
As before, $\delta: X\to \bP^1 \times \bP^2$ is a double cover ramified along a smooth divisor of bi-degree $(2, 2)$.

\begin{proposition}
\label{proposition:example_involution}
    Let $\tau \in \mathrm{Aut}(X)$ be an involution such that
    \begin{itemize}
        \item the induced action on $\bP^1$ via $\pi_1: X \to \bP^1$ is trivial,
        \item the induced action on $\bP^2$ via $\pi_2 : X \to \bP^2$ is non-trivial, and
        \item $\widetilde{\Delta}^\tau = \emptyset$.
    \end{itemize}
    Then $X$ is not $\langle \tau \rangle$-projectively linearizable.
\end{proposition}

\begin{proof}
    Since cyclic actions on projective space are linearizable, we may assume that the induced action on $\mathbb P^2$ is given by
    \[
    \bar{\tau}: [x:y:z] \mapsto [x:y:-z].
    \]
    Then the discriminant curve has the form  
    \[
    z^4 + z^2Q_2(x, y) + Q_4(x, y),
    \]
    where $Q_2, Q_4\in \bC[x,y]$ are homogeneous of degree $2$ and $4$, respectively. Since we assume that $\Delta$ is smooth, $Q_4 = 0$ has $4$ distinct roots.
    
    Suppose that $X$ is $\langle \tau \rangle$-projectively linearizable. By Theorem~\ref{theorem:IJTobstructions}, there exists a smooth connected projective curve $C$ of genus $2$ with $\fS_2$-action such that $(\mathrm{\mathbf{CH}}^2_{X/\mathbb C})^{\gamma_1 + \gamma_2}$ is equivariantly isomorphic to $\mathrm{Pic}^m(C)$, for some $m$.
    The quotient of $C$ by a nontrivial involution is either 
   an elliptic curve, or $\mathbb P^1$.
    In either case, $C$ has a fixed point, and thus $\Pic^m(C)$ is equivariantly isomorphic to $\Pic^0(C)$, with a fixed point. 

    We conclude that $(\mathrm{\mathbf{CH}}^2_{X/\mathbb C})^{\gamma_1 + \gamma_2}$ has a fixed point. 
    Since 
    $$\cF_{1, 1}(X) \to (\mathrm{\mathbf{CH}}^2_{X/\mathbb C})^{\gamma_1 + \gamma_2}
    $$ 
    is a $\bP^1$-fibration, $\cF_{1,1}(X)$ also has a fixed point.
    Let $f :R \to X$ be a stable map which is $\langle \tau \rangle$-stable. If $R$ is reducible, then $f(R)$ consists of two lines in fibers of $\pi_1$ and $\pi_2$ respectively. However, $\tau$ never preserves a line in fibers of $\pi_2$ because $\widetilde{\Delta}^\tau = \emptyset$.  
    Thus $R$ is irreducible and $R_1 = f(R)$ is $\langle \tau \rangle$-stable. If $R_1$ is not fixed, then $R_1$ must be a bisection of $\pi_1$ because the induced action on $\mathbb P^1$ is trivial. This contradicts the fact that the class of $R_1$ is $\gamma_1 + \gamma_2$. We conclude that $R_1$ is fixed. Let $\ell$ be the image of $R_1$ via $\pi_2 : X \to \mathbb P^2$. Then $\ell$ is a line, fixed by the group action, i.e., 
    $\ell$ is defined by $z = 0$.
    Let $S$ be the pullback of $\ell$ via $\pi_2 : X \to \mathbb P^2$.
    Then $S$ is a smooth del Pezzo surface of degree $4$ because $Q_4 = 0$ has distinct roots.
    
    Consider the double cover $S \to \mathbb P^1 \times \ell$ and let $R'$ be 
   the image of $R_1$ via this map. Then $R_1 \to R'$ is birational. This means that the pullback of $R'$ via $S \to \mathbb P^1 \times \ell$ consists of two curves $R_1$ and $R_2$ and the involution is swapping those curves. This is a contradiction.
\end{proof}

\begin{example}
\label{example:involutions_conicbundle}
    Consider the double cover:
$$
X: w^2 = t_0^2q_0 + 2t_0t_1q_{01} + t_1^2q_1,
$$
with 
\begin{align*}
  q_0   &:=  z^2 + a_{0,0} x^2 + a_{0, 1} xy + a_{0, 2} y^2, \\
  q_{01}&:=a_{01, 0} x^2+a_{01, 1} xy + a_{01, 2} y^2, \\
  q_1&:=z^2 + a_{1, 0} x^2 + a_{1, 1} xy + a_{1, 2} y^2,
\end{align*}
    where $a_{i, j} \in \mathbb C$.
    The discriminant curve $\Delta$ is given by
    \[
    z^4 + z^2Q_2(x, y) + Q_4(x, y) = 0,
    \]
    where $Q_2,Q_4\in \bC[x,y]$ are homogeneous polynomials of degree $2$ and $4$, respectively. We assume that $\Delta$ is smooth and 
    consider the involutions
\begin{align*} 
    \mu : ([t_0:t_1], [x:y:z], w) & \mapsto ([t_0:t_1], [x:y:-z], w), \\
    \sigma : ([t_0:t_1], [x:y:z], w) & \mapsto ([t_0:t_1], [x:y:z], -w).   
\end{align*} 
The involution $\mu$   
is linearizable by Proposition~\ref{proposition:prokhorov2018}, since $\widetilde{\Delta}$ has a $\langle \mu \rangle$-fixed point. The action of $\tau := \sigma \mu$  is not projectively linearizable, by Proposition~\ref{proposition:example_involution}.
\end{example}

\begin{remark}
The Burnside formalism of \cite{BnG} applies in Example~\ref{example:involutions_conicbundle} to $G\simeq \fS_2^2=\langle\sigma, \tau\rangle$: 
the $\sigma$-action fixes a del Pezzo surface of degree 2, and the residual action fixes a curve of genus 3. The formalism of incompressible symbols of, e.g., \cite[Section 3.6]{TYZ23} shows that the $G$-action on $X$ is not projectively linearizable.
\end{remark}

    \begin{example}
    \label{example:C_2timesC_2}
Consider the double cover:
$$
X: w^2 = t_0^2q_0 + 2t_0t_1q_{01} + t_1^2q_1,
$$
with 
\begin{align*}
  q_0   &:=  a_{0, 0} x^2 + a_{0, 1} y^2+ a_{0,2} z^2, \\
  q_{01}&:= a_{01,0} x^2+a_{01, 1} y^2 + a_{01, 2} z^2,\\
  q_1&:=a_{1, 0} x^2 + a_{1, 1} y^2 + a_{1, 2} z^2,
\end{align*}
where $a_{i, j} \in \mathbb C$.
In this case, the discriminant $\Delta\subset \bP^2$ is defined by an equation of the form
    \[
    Q(x^2, y^2, z^2) = 0.
    \]
We assume that $\Delta$ is a smooth quartic and consider the involutions
\begin{align*}
&\tau_1 : ([t_0:t_1], [x:y:z], w) \mapsto ([t_0:t_1], [-x:y:z], w)\\
&\tau_2 : ([t_0:t_1], [x:y:z], w) \mapsto ([t_0:t_1], [x:-y:z], w)\\
&\tau : ([t_0:t_1], [x:y:z], w) \mapsto ([t_0:t_1], [x:y:-z], -w).
\end{align*}
Then $\tau_1\tau_2 = \tau$ and $G:=\langle \tau_1,\tau_2\rangle \simeq \fS_2^2$
does not contain the involution $\sigma$ from Example~\ref{example:involutions_conicbundle}. 
For $i = 1, 2$, we have $\widetilde{\Delta}^{\tau_i} \neq \emptyset$, so that $X$ is $\langle \tau_i \rangle$-projectively linearizable.
However, by Proposition~\ref{proposition:example_involution}, $X$ is not $\langle \tau \rangle$-projectively linearizable, and thus not $G$-projectively linearizable.
    \end{example}

      \begin{remark}
        In Examples~\ref{example:involutions_conicbundle} and \ref{example:C_2timesC_2}, (projective) linearizability is equivalent to the vanishing of equivariant IJT obstructions by Theorem~\ref{theorem:rationalitycriterion}. It follows from Corollary~\ref{thm:distinguish} that equivariant birational classes of conic bundles $X, X'$ with involutions $\tau$, respectively, $\tau'$, are different as soon as $C, C'$ are not equivariantly isomorphic.
        In these examples, there are no non-uniruled divisors in the fixed locus of the involutions; previously known obstructions do not apply to these conic bundles, see \cite{pro-inv}. 
    \end{remark}

The following example shows that $\widetilde{\Delta}^G \neq \emptyset$ is not a necessary condition for linearizability.

\begin{example}
\label{example:conic_sigma}
Here we prove that $X$ is always $\langle \sigma \rangle$-linearizable, where $\sigma$ is the unique  nontrivial involution acting trivially on $\mathbb P^1 \times \mathbb P^2$.
This follows from the triviality of $\widetilde{P}^{(1)}$, as a $P$-torsor by Proposition~\ref{proposition:11}.

Let $\ell$ be a bitangent line to $\Delta$ such that $\ell\cap \Delta$ consists of two points. Then $S_\ell$ contains two $A_1$-singularities, and there is a curve $R$ of class $\gamma_1 + \gamma_2$ mapping to $\ell$ and passing through these two singularities. 
Indeed, let $\widetilde{S}_\ell \to S_\ell$ be the minimal resolution. A singular fiber of $\pi_1|_{S_\ell} : S_\ell \to \bP^1$ at $t_0 = 0$ contains a line $\ell'$ of class $\gamma_1$, and such a line passes through one singular point but not the other. We also have $-K_{S_\ell}.\ell' = 1$. The fibration $\widetilde{S}_\ell \to \ell$ comes with two singular fibers, both consisting of two $(-1)$-curves connected by one $(-2)$-curve. If we contract all $(-1)$-curves of these singular fibers, we get a weak del Pezzo surface of degree $8$.
    The strict transform of $\ell'$ is a conic of self-intersection $0$. This implies that this weak del Pezzo surface is a smooth quadric surface. The existence of $R$ follows.
Then the existence of $R$ induces a fixed point on $\widetilde{P}^{(1)}$. 

When $\ell\cap \Delta$ consists of one point, $S_\ell$ contains a unique singularity of type $A_3$. The minimal resolution $\widehat{S}_\ell$ of $S_\ell$ contains a chain of three $(-2)$-curves, with two $(-1)$-curves attached to the edges. Then there is a curve of class $\gamma_1 + \gamma_2$ mapping to $\ell$ whose strict transform on $\widehat{S}_\ell$ meets the central $(-2)$-curve of the chain. This induces a fixed point on $\widetilde{P}^{(1)}$. 
\end{example}

\section{Quadric surface bundles}
\label{sect:quad}

Let $\cX$ be a smooth projective threefold over $\bC$ admitting a quadric surface bundle structure
$$
\pi : \mathcal X \to \bP^1.
$$
Such threefolds are rational: by Tsen's theorem, there exists a section of $\pi$, projecting from it gives birationality to a relative $\bP^2$-bundle, thus birationality $X\sim \bP^3$. However, this may not work equivariantly. 

Here, we explore equivariant birational geometry of such quadric surface bundles, establish a criterion for linearizability of actions of cyclic groups, and
provide examples of nonlinearizable actions. 

\subsection*{Geometry} 
We follow \cite{HT12}. The intermediate Jacobian is given by
$$
\IJ(\cX)=\Pic^0(C),
$$
where $C\to \bP^1$ is the discriminant double cover, a
hyperelliptic curve ramified in the discriminant $\mathfrak d\subset \bP^1$ of $\pi$ (reduced, since $\cX$ is smooth, by assumption). Let
$|\mathfrak d|$ be the degree of the discriminant divisor,  
\[
g(C) = \frac{|\mathfrak d|}{2}-1
\]
the genus of $C$, and $\omega_\pi$ the relative dualizing sheaf of $\pi$.
The height of a section $s : \bP^1 \to \mathcal X$ is defined as
\[
h_{\omega_\pi^{-1}}(s) := \deg (s^*\omega_\pi^{-1}),
\]
it is equal to the degree of the normal bundle $N_s$.
The space of sections of $\pi$ of height $h$ is denoted by
\[
\mathrm{Sec}(\mathcal X/\bP^1, h),
\]
and its closure in the coarse moduli space of stable maps of genus $0$ by
\[
\overline{\mathrm{Sec}}(\mathcal X/\bP^1, h).
\]
The expected dimension of this moduli space is $h + 2$.
Let $\cF_1(\mathcal X/\bP^1)$ the relative Fano variety of lines, with the Stein factorization
\[
\cF_1(\mathcal X/\bP^1) \to C \to \bP^1.
\]

\begin{lemma}[{\cite[Proposition 2]{HT12}}]
\label{lemma:HT12Proposition2}
Assume that $C \to \bP^1$ is non-split, i.e., $C$ is irreducible.
Fix a class $\gamma$ of sections of $\pi$ of  sufficiently large height. Then
    the Abel-Jacobi map
    \[
    \mathrm{AJ} : \mathrm{Sec}(\mathcal X/\bP^1, \gamma) \to (\mathrm{\mathbf{CH}}^2_{X/\mathbb C})^\gamma
    \]
    is an open immersion composed with a projective bundle over $(\mathrm{\mathbf{CH}}^2_{X/\mathbb C})^\gamma$.
\end{lemma}

\begin{proof}
    We include a proof,
    following \cite{HT12}.
    Since the Brauer group of $\bC(C)$ is trivial, $\mathcal F_1(\mathcal X/\bP^1) \to C$ is the projectivization $\mathbb P_C(V)$ of a rank $2$ vector bundle $V$ over $C$. Let $\mathcal O_{\mathbb P_C(V)}(1)$ be the relative polarization of $\mathbb P_C(V)$. There is a bijection between sections of $\mathcal X/\bP^1$ and sections of $\mathcal F_1(\mathcal X/\bP^1)/C$, as in \cite[Section 3]{HT12},
    with a natural identification 
    \[
    \mathrm{Sec}(\mathcal X/\bP^1, h) = \mathrm{Sec}(\mathcal F_1(\mathcal X/\bP^1)/C, d),
    \]
    where $d = \deg(s^*\mathcal O_{\mathbb P_C(V)}(1))$ and $s : C \to \mathcal F_1(\mathcal X/\bP^1)$ is a section. 
    By \cite[(3.2)]{HT12},
    \[
     h = 2d + \deg (V) - \frac{|\mathfrak d|}{2}.
    \]
    We have a natural morphism
    \[
    \mathrm{Sec}(\mathcal F_1(\mathcal X/\bP^1)/C, d) \to \Pic^d(C), \quad  s \mapsto s^*\mathcal O_{\mathbb P_C(V)}(1).
    \]
    Let $\mathcal L$ be the universal line bundle on $C \times \Pic^d(C)$.
    When $d$ is sufficiently large, $(\rho_2)_*(\mathcal L\otimes \rho_1^*V)$ is locally free by \cite[III Theorem 12.11]{Hartshorne}; here $\rho_i$ are the projections onto the factors in $C \times \Pic^d(C)$.  
    By \cite[Proposition 2]{HT12}, the moduli space
    \[
    \mathrm{Sec}(\mathcal F_1(\mathcal X/\bP^1)/C, d)
    \]
    is a Zariski open in the projectivization $\mathbb P((\rho_2)_*(\mathcal L\otimes \rho_1^*V))$ over $\Pic^d(C)$. 

    To show that this morphism coincides with the Abel-Jacobi map, we follow \cite[Example 10.4]{LT24}: 
    after replacing $V$ by $V\otimes L$ with $L$ being a line bundle on $C$, we may assume that $\mathcal O_{\mathbb P_C(V)}(1)$ admits a section $C'$; we let $B'$
    be a section of $\pi$ corresponding to $C'$. We denote the blow up of $\mathcal X$ along $B'$ by $\phi : \mathcal X' \to \mathcal X$. Let $E \subset \mathcal X$ be the divisor swept out by lines of $\pi$ meeting at $B'$, we denote its strict transform by $E' \subset \mathcal X'$. This is a $\mathbb P^1$-bundle over $C'$ and we denote this projection by $\eta : E' \to C'$.
    Consider the homomorphism
    \[
    \rH_3(\mathcal X, \mathbb Z) \stackrel{\phi^*}{\lra} \rH_3(\mathcal X', \bZ) \stackrel{\cap E'}{\lra} \rH_1(E', \mathbb Z) \stackrel{\eta_*}{\lra} \rH_1(C', \mathbb Z),
    \]
    where $\phi^*$ is the homology pullback defined using the Poinc\'are duality twice. It respects the Hodge structures, giving a homomorphism
    \[
    \rH^{2, 1}(\mathcal X)^\vee \to \rH^{2, 1}(\mathcal X')^\vee \to \rH^{1, 0}(E')^\vee \to \rH^{1, 0}(C')^\vee.
    \]
    Fix homologous sections $B_0, B_1$ of $\pi$ and denote the strict transforms of $B_i$ via $\beta$ by $B_i'$.
    There is a real $3$-cycle $\alpha$ such that $\partial \alpha = B_1 - B_0$ and  $\partial \phi^*\alpha = B_1' - B_0'$. The Abel-Jacobi map sends $B_1$ to 
    $$
    \int_\alpha  \cdot \in \rH^{2, 1}(\mathcal X)^\vee/\rH_3(\mathcal X, \mathbb Z).
    $$
    This is mapped by the above homomorphism to $\int_{\phi^*\alpha \cap E'} \cdot$ and $\int_{\eta_*(\phi^*\alpha \cap E')}\cdot$. The boundary of $\eta_*(\phi^*\alpha \cap E)$ is given by $$
    \eta(B_1'\cap E') - \eta(B'_0\cap E').
    $$
    This gives rise a morphism
    \[
    \mathrm{Sec}(\mathcal X/\bP^1, \gamma) \to \mathrm{IJ}(\mathcal X) \to \mathrm{Jac}(C) \cong \mathrm{Pic}^0(C).
    \]
    This composition is compatible with the construction of \cite[Proposition 2]{HT12} described above. On the other hand, since the fibers are open subsets of projective spaces, the Abel-Jacobi map must factor as
    \[
    \mathrm{Sec}(\mathcal X/\bP^1, \gamma) \to \mathrm{Pic}^0(C) \to \mathrm{IJ}(\mathcal X).
    \]
\end{proof}

\subsection*{Automorphisms}

Let $G\subseteq \Aut(\cX)$ be a finite group. Assume 
that $\pi$ is equivariant for $G$.
Easy examples of projectively nonlinearizable actions on $\cX=\bP^1\times Q$, where $Q$ is a smooth quadric surface, are:
\begin{itemize}
\item generically free action of $\fS_2^2$ on $\bP^1$ 
and any action on $Q$ -- this is not linearizable since there are no fixed points on $\cX$; 
\item $\fS_2$ acting on $\bP^1$ and $\mathfrak D_6$ acting on $Q$ -- via the Burnside formalism of \cite{BnG}, see \cite[Theorem 5.1]{CTZ-segre}. 
\end{itemize}

\subsection*{The main example}
Consider 
$$
\cX\subset \bP^1\times \bP^3, 
$$
given by the vanishing of a sufficiently general form $F$ of bi-degree $(n,2)$. Projection to the first factor $\cX\to \bP^1$ is a quadric surface bundle of relative Picard rank $1$,  
splitting over the discriminant double cover $C\to \bP^1$ and ramified over 
the discriminant divisor $\mathfrak d$ of degree $|\mathfrak d|=4n$; 
so that $g(C)=2n-1$.

\subsection*{Intermediate Jacobians and torsors}
Let $\pi : \mathcal X \to \bP^1$ be a smooth quadric surface bundle with a regular, compatible, action by a finite group $G$. This induces a $G$-action on $\cF_1(\mathcal X/\bP^1)$ so that
\[
\cF_1(\mathcal X/\bP^1) \to \bP^1,
\]
and the Stein factorization
\[
\cF_1(\mathcal X/\bP^1) \to C \to \bP^1
\]
are $G$-equivariant. 
Let $\gamma_\ell$ be the class of lines in fibers of $\pi$. 

\begin{lemma}
\label{lemma:IJ_quadric}
    Assume that $C \to \bP^1$ is non-split.
    Then $\gamma_\ell$ is $G$-invariant and
    we have a $G$-equivariant isomorphism
    \[
    (\mathrm{\mathbf{CH}}^2_{\mathcal X/\mathbb C})^{\gamma_\ell} \cong \Pic^1(C).
    \]
    Moreover, we have a $G$-equivariant isomorphism
    \[
    \IJ(\mathcal X) \cong (\mathrm{\mathbf{CH}}^2_{\mathcal X/\mathbb C})^{0} \cong \Pic^0(C),
    \]
    as $G$-equivariant principally polarized abelian varieties.
\end{lemma}

\begin{proof}
    Note that 
    $
    \cF_1(\mathcal X/\bP^1) \to C
    $    is a $\mathbb P^1$-bundle over $C$.
    The Abel-Jacobi map
    \[
    \mathrm{AJ} : \cF_1(\mathcal X/\bP^1) \to (\mathrm{\mathbf{CH}}^2_{\mathcal X/\mathbb C})^{\gamma_\ell},
    \]
    factors as
    \[
    \cF_1(\mathcal X/\bP^1) \to C \to (\mathrm{\mathbf{CH}}^2_{\mathcal X/\mathbb C})^{\gamma_\ell}
    \]
     and induces a $G$-equivariant isomorphism
    \[
    \Pic^1(C) \cong (\mathrm{\mathbf{CH}}^2_{\mathcal X/\mathbb C})^{\gamma_\ell},
    \]
as in the proof of Lemma~\ref{lemma:HT12Proposition2}.

The claim concerning principal polarizations follows from the fact that, in the notation of the proof of Lemma~\ref{lemma:HT12Proposition2}, $\mathcal X'$ is the blowup of a $\bP^2$-bundle over $\bP^1$ along a curve $C'$. Indeed, one can use the blowup formula to verify compatibility of cup products.
\end{proof}

\subsection*{Linearization criterion}
 Over non-closed fields, if $\pi : X \to \mathbb P^1$ admits a section defined over $k$, then $X$ is $k$-rational. Here is an analogous statement in the equivariant setting:
\begin{lemma}
\label{lemma:rationality_quadric}
Let $G$ be a cyclic group and 
$\pi : \mathcal X \to \bP^1$ a $G$-equivariant quadric surface bundle with a $G$-invariant section. Then the $G$-action on $\cX$ is linearizable.
\end{lemma}

\begin{proof}
    Projecting from a $G$-invariant section gives an equivariant birationality of $\cX$ to a $\mathbb P^2$-bundle over $\bP^1$. Since $G$ is cyclic, the generic fiber $\mathbb P^2_\eta$ admits a $G$-invariant line. Such a $\mathbb P^2_\eta$ is linearizable, and we obtain a $G$-equivariant birationality $\cX\sim_G \bP^1\times \bP^2$, and thus linearizability. 
    \end{proof}

\begin{theorem}
\label{theorem:quadric_noneffective}
Let $\pi:\mathcal X\to \bP^1$
be a smooth quadric surface bundle of 
relative Picard rank $1$ and $|\mathfrak d| \geq 6$. 
Let $G\subseteq \Aut(\cX)$ be a cyclic group such that $\pi$ is $G$-equivariant. 
Then the $G$-action on $\cX$ is
linearizable if and only if there exists a $G$-invariant section of $\pi$.
\end{theorem}

\begin{proof}
If $\pi$ admits a $G$-invariant section, 
Lemma~\ref{lemma:rationality_quadric} applies. 
Conversely, assume that the $G$-action on $\cX$ is linearizable.
    Fix a class $\gamma$ of sections. Since $\cX$ has Picard rank $2$ and $G$ respects $\pi$, $\gamma$ is a $G$-invariant class.
    By Theorem~\ref{theorem:IJTobstructions}, we have
    \[
    (\mathrm{\mathbf{CH}}^2_{\mathcal X/\mathbb C})^\gamma \cong  \Pic^e(C),
    \]
    for some $e$. By Lemma~\ref{lemma:IJ_quadric}, we have
    \[
    [(\mathrm{\mathbf{CH}}^2_{\mathcal X/\mathbb C})^{\gamma + k\gamma_\ell}] = [e + k][\Pic^1(C)].
    \]
    In particular, whenever $e + k \equiv 0 \pmod{ 2g(C)-2}$, this torsor is trivial.

    Let $k$ be sufficiently large so that the height of $\gamma + k\gamma_\ell$ is large and $e + k \equiv 0 \pmod{2g(C)-2}$.
    By Lemma~\ref{lemma:HT12Proposition2}, the Abel-Jacobi map
    \[
    \mathrm{AJ} : \mathrm{Sec}(\mathcal X/\bP^1, \gamma + k\gamma_\ell) \to (\mathrm{\mathbf{CH}}^2_{\mathcal X/\mathbb C})^{\gamma + k\gamma_\ell}
    \]
    is a composition of an open immersion followed by a projective bundle 
    $$
    \mathcal P \to (\mathrm{\mathbf{CH}}^2_{\mathcal X/\mathbb C})^{\gamma + k\gamma_\ell}.
    $$
    Since $(\mathrm{\mathbf{CH}}^2_{\mathcal X/\mathbb C})^{\gamma + k\gamma_\ell}$ has a fixed point, and $G$ is cyclic, $\mathcal P$ has a fixed point as well.
    This implies that $\overline{\mathrm{Sec}}(\mathcal X/\bP^1, \gamma + k\gamma_\ell)$ has a fixed point, proving that there is a $G$-invariant section of $\pi$. 
\end{proof}

\subsection*{Effective results}

Here we establish an effective version of Theorem~\ref{theorem:quadric_noneffective}.
First we recall:

\begin{proposition}[{\cite[Proposition 15]{HT12}}]
\label{proposition:HT12Proposition14}
    Let $\pi : \mathcal X \to \bP^1$ be a smooth quadric surface bundle with discriminant of degree $|\mathfrak d|$ such  that
    \begin{itemize}
        \item the discriminant cover $C \to \bP^1$ is non-trivial, and
        \item for any section $s : \bP^1 \to \mathcal X$, we have
        \[
         h_{\omega_\pi^{-1}}(s) \geq -\frac{|\mathfrak d|}{2}.
        \]
    \end{itemize}
    Then $\mathcal F_1(\mathcal X/\bP^1) \to C$ is the projectivization of a semi-stable bundle $V$ over $C$.
\end{proposition}

After rescaling, we may, and will, assume that $\deg (V) = 0$ or $1$.

\begin{proposition}
\label{proposition:effectiveHT12proposition2}
With the assumptions of Proposition~\ref{proposition:HT12Proposition14}, 
consider  
$$
\bP_C(V)=\mathcal F_1(\mathcal X/\bP^1) \to C
$$ 
Suppose that $|\mathfrak d| \geq 6$ and $\gamma$ is a class of sections of height $h$. Assume that
$$
h \geq \frac{3}{2}|\mathfrak d| -4 + 4\deg (V).
$$
Then the statement of Lemma~\ref{lemma:HT12Proposition2} is valid. 
\end{proposition}

\begin{proof}
    We use notation from Lemma~\ref{lemma:HT12Proposition2}.
    By Proposition~\ref{proposition:HT12Proposition14}, $V$ is semi-stable. 
    By Riemann-Roch, there exists a section 
    $$
    s_0 : C \to \mathbb P_C(V)
    $$ 
    with  $$
    d_0 = \deg(s^*\xi)  \leq g(C) = \frac{|\mathfrak d|}{2}-1,
    $$ 
    where $\xi = \mathcal O_{\mathbb P_C(V)}(1)$.
    The numerical class of the canonical divisor of $\mathbb P_C(V)$ is given by
    \[
    K_{\mathbb P_C(V)} = -2\xi + (2g(C)-2 + \deg (V))\lambda   ,
    \]
    where $\lambda$ is the class of fibers of $\mathbb P_C(V) \to C$. Thus
    the degree of $N_{s_0}$ is 
    \[
    2d_0 - \deg (V).
    \]
    We then glue 
    $$
    b \geq 2g(C) -2d_0 +\deg (V) = |\mathfrak d|-2 -2d_0 + \deg (V)
    $$ 
    fibers of $\mathbb P_C(V) \to C$ to $s_0(C)$.
    By the assumption of Proposition~\ref{proposition:HT12Proposition14} and \cite[(3.2)]{HT12}, 
    $$
    d_0 = \frac{h_0}{2} - \frac{\deg (V)}{2} + \frac{|\mathfrak d|}{4} \geq - \frac{\deg (V)}{2}.
    $$
    Using \cite[Lemma 2.6]{GHS} and smoothing arguments,
    the existence of sections is guaranteed as soon as 
    \[
    d \geq |\mathfrak d|  -2+ \frac{3}{2}\deg (V),
    \]
    which in turn follows when
    \[
    h \geq \frac{3}{2}|\mathfrak d| -4 + 4\deg (V).
    \]
    Assuming 
    $
    d \geq 2g(C),
    $ 
    it follows from Proposition~\ref{proposition:HT12Proposition14} and \cite[Corollary 2.8]{LRT23} that
    $
    V \otimes s^*\mathcal O_{\mathbb P_C(V)}(1)
    $
    is globally generated and 
    $$
    \rH^1(C, V \otimes s^*\mathcal O_{\mathbb P_C(V)}(1)) = 0,
    $$ 
    where $s : C \to \mathbb P_C(V)$ is a section such that $\deg( s^*\xi) = d$.
    This implies that $(\rho_2)_*(\mathcal L\otimes \rho_1^*V)$ is locally free. (Note that this is the only requirement for the proof of Lemma~\ref{lemma:HT12Proposition2}.)
    The condition 
     $
    d \geq 2g(C)
    $ 
    translates to 
    \[
    h \geq 4g(C) + \deg (V) -\frac{|\mathfrak d|}{2} = \frac{3}{2}|\mathfrak d| + \deg (V) -4.
    \]
\end{proof}

Here is an effective version of Theorem~\ref{theorem:quadric_noneffective}:

\begin{theorem}
\label{theorem:quadric_effective}
Let $G$ be a cyclic group and  
$$
\pi:\mathcal X\to \bP^1
$$
a $G$-equivariant quadric surface bundle of 
relative Picard rank $1$ and 
$|\mathfrak d| \geq 6$. With the assumptions of Proposition~\ref{proposition:effectiveHT12proposition2}, the $G$-action on 
$\cX$ is linearizable if and only if there exists a $G$-invariant section of height 
$$
\leq \frac{7}{2}|\mathfrak d| + 4\deg (V)  - 14.
$$
\end{theorem}

\begin{proof}
    Assume that $\mathcal X$ is linearizable.
    In the proof of Theorem~\ref{theorem:quadric_noneffective}, the height $h(\gamma) + 2k$ of $\gamma + k\gamma_\ell$ needs to be $\ge \frac{3}{2}|\mathfrak d| -4+ 4\deg (V)$ and we also must have $e + k \equiv 0 \mod (|\mathfrak d|-4)$.
    One can find $k$ such that 
    \[
    \frac{3}{2}|\mathfrak d| -4+ 4\deg (V) \leq h(\gamma) + 2k\leq  \frac{7}{2}|\mathfrak d| + 4\deg (V) - 14,
    \]
    and $e + k \equiv 0 \mod (|\mathfrak d|-4)$.
    Thus our assertion follows.
\end{proof}

\begin{example}
\label{example:(n, 2)_effecitve}
Let $\mathcal X \subset \bP^1 \times \mathbb P^3$ be a smooth divisor of type $(n, 2)$, where $n \geq 2$; with $\pi : \mathcal X \to \bP^1$ a quadric surface bundle of relative Picard rank $1$ with square-free discriminant and $|\mathfrak d| = 4n \geq 8$.
By \cite[Section 4, Case 1]{HT12}, $\deg (V) \equiv n \pmod{2}$.

    Any section $s : \bP^1 \to \mathcal X$ has height of the form $2a-n$, with $a = \deg(s^*H) \geq 0$, where $H$ is the pullback of the hyperplane class from $\mathbb P^3$. Indeed, this follows from the adjunction formula and the fact that the relative dualizing divisor is $K_{\mathcal X} - K_{\bP^1}$. Since 
    \[
    -n \geq - \frac{|\mathfrak d|}{2} = -2n,
    \]
    the assumptions of Proposition~\ref{proposition:HT12Proposition14} are verified. By Theorem~\ref{theorem:quadric_effective}, 
    a cyclic action on $\pi: \mathcal X\to \bP^1$ is projectively linearizable if and only if there exists a $G$-invariant section of height 
    $
    \leq 14n + 4 \deg (V) -14.
    $
\end{example}

\subsection*{Examples of actions}
We consider $\mathcal X \subset \mathbb P^1 \times \mathbb P^3$, with $n \geq 2$, defined by the 
vanishing of a form $F$ of bi-degree $(n, 2)$. Let 
$$
\pi:\cX\to \bP^1
$$ 
be the associated quadric surface bundle.

\begin{example}
\label{example:fixedcurves}
Let 
\[
F:=\sum_{i = 0}^n t_0^{i}t_1^{n-i}(f_i(x, y) + g_i(z, w)),
\]
where $f_i, g_i$ are general binary quadratic forms, 
so that $\mathfrak d$ is reduced, and thus $\mathcal X$ smooth.
The involution
\[
\sigma : ([t_0:t_1], [x:y:z:w]) \mapsto ([t_0:t_1], [-x:-y:z:w])
\]
acts trivially on $\mathbb P^1$, and any invariant section of $\pi$ is pointwise fixed. However, the fixed locus of $\sigma$ consists of two nonrational curves, thus there are no invariant sections.
By Theorem~\ref{theorem:quadric_noneffective}, 
the $\langle\sigma\rangle$-action on 
$\mathcal X$ is not (projectively) linearizable.
\end{example}

\begin{example}
\label{example:fixedquadric}
Let 
    \[
    F:=t_0^2f(x, y, z, w) + t_1^2g(x, y, z, w),
    \]
    where $f, g$ are general quadratic forms in $4$ variables.
Consider
\[
\sigma : ([t_0:t_1], [x:y:z:w]) \mapsto ([-t_0:t_1], [x:y:z:w]).
\]
Note that $\sigma$ is the covering involution of $\pi_2: \mathcal X \to \bP^3$. If we had a $\langle\sigma\rangle$-invariant section $R \subset \mathcal X$, then $\pi_2|_{R} : R \to \pi_2(R) = R'$ would be of degree $2$. Since $R$ is smooth, $R'$ must be smooth as well, and the branch divisor of $R \to R'$ reduced, showing that $R$ cannot be rational, contradiction. Thus,
the $\langle\sigma\rangle$-action is not projectively linearizable.
\end{example}

\begin{remark} 
Examples~\ref{example:fixedcurves} and \ref{example:fixedquadric} produce moduli of conjugacy classes of involutions in $\mathrm{Cr}_3(\bC)$, by  varying the underlying curve $C$ in the moduli space and applying Corollary~\ref{thm:distinguish}. In these examples, the non-uniruled divisorial component of the fixed locus of the involution $\sigma$ is empty, and the previously known obstructions do not distinguish these equivariant birational classes.
\end{remark}

\section{Fano threefolds}
\label{sect:fano}

Let $X$ be a smooth projective Fano threefold of Picard rank 1 
and $\cF_d(X)$ the variety of rational curves of anticanonical degree $d$. 
Over nonclosed fields, rationality of geometrically rational varieties of this type has been recently settled in \cite{HT-odd}, \cite{HT-cycle}, \cite{BW19}, \cite{HT-18}, \cite{KuzP}. There are three types of such Fano varieties:
\begin{itemize}
    \item rational over any field: the quintic del Pezzo threefold $V_5\subset \bP^4$, 
    \item rational if and and only if there are rational points: forms of $\bP^3$, the smooth quadric $Q\subset \bP^4$, $X_{12}\subset \bP^8$ and $X_{22}\subset \bP^{13}$, 
    \item rational if and only if there are rational points and rational curves, defined over the ground field, in prescribed degrees: 
        \begin{itemize}
            \item $X_{2,2}\subset \bP^5$, intersection of two quadrics -- lines,
            \item $X_{16}\subset \bP^{10}$ -- twisted cubics,
            \item $X_{18}\subset \bP^{13}$ -- conics. 
        \end{itemize}
\end{itemize}
In \cite{KuzP}, the intermediate Jacobian torsor obstructions of \cite{HT-quad}, \cite{HT-cycle}, and \cite{BW19} are reinterpreted in the framework of derived categories. 

The equivariant situation is markedly different. First of all, the quintic del Pezzo threefold  $V_5$ is {\em not} projectively linearizable for the action of  $\fA_5\subset \Aut(V_5)=\PGL_2(\bC)$, by the main result of \cite{CS}. In fact, $V_5$ is birationally rigid for this action. 
Projective linearizability of quadrics is an open problem; linearizability is unknown even for some actions of $\fS_3$. In \cite[Section 5]{BBT} there is an example of a nonlinearizable action by the Frobenius group $\mathfrak F_8$ on $X_{12}$ without any obstructions in the derived category. 

The following theorem is a partial equivariant analog:

\begin{theorem} 
\label{thm:fanolinear}
Let $X$ be a smooth projective Fano threefold 
of the following types:
\begin{itemize}
    \item[(1)] an intersection of two quadrics $X_{2,2}\subset\bP^5$, 
    \item[(2)] a prime Fano threefold of genus 9, $X_{16}\subset \bP^{10}$,
    \item[(3)] a prime Fano threefold of genus 10, $X_{18}\subset \bP^{11}$. 
\end{itemize}
Let $G$ be a finite group acting generically freely and regularly on $X$. Assume that the $G$-action is projectively linearizable. Then 
\begin{equation} 
\label{eqn:cond}
\cF_d(X)^G\neq \emptyset,
\end{equation}
for $d=2,3,$ or $2$, respectively. 

In Cases (1) and (2), this condition is also sufficient for projective linearizability. 
\end{theorem}

We treat the cases separately.

\subsection*{Intersections of two quadrics}

This case was settled in \cite{HT-odd}, via reduction to nonclosed fields. We supply an argument based on the equivariant intermediate Jacobian torsor formalism of Section~\ref{sect:jacob}. After the completion of this paper, F. Scavia sent us his unpublished note \cite{Scavia} with a similar proof. 

The variety of lines $\cF_2(X_{2,2})$ is a torsor under $\IJ(X)=\Pic^0(C)$. 
When if $\cF_2(X_{2, 2})^G \neq \emptyset$, there exists a $G$-stable line, which yields a standard $G$-equivariant birationality with $\bP^3$. 
To prove the converse, we need: 

\begin{lemma}
\label{lemma:X22}
    Let $X_{2, 2}$ be a smooth complete intersection of two quadrics in $\mathbb P^5$ with a regular $G$-action on $X_{2, 2}$. Then the Abel-Jacobi map defines a $G$-equivariant isomorphism
    \[
    \mathrm{AJ}: \cF_2(X_{2, 2}) \cong (\mathrm{\mathbf{CH}}^2_{X_{2, 2}/\mathbb C})^\ell,
    \]
    where $\ell$ is the class of lines.
    Moreover, as a $G$-equivariant torsor, we have
    \[
    2[\cF_2(X_{2, 2})] = [\mathrm{Pic}^1(C)],
    \]
    where $C$ is a smooth projective curve of genus $2$ with a regular $G$-action.
\end{lemma}

\begin{proof}
    The non-trivial part of the first statement is the $G$-equivariant property, established in Theorem~\ref{thm:tors}.
For the second statement, we need to show that we have a $G$-equivariant isomorphism
    \[
    (\mathrm{\mathbf{CH}}^2_{X_{2, 2}/\mathbb C})^{2\ell} \cong \Pic^1(C),
    \]
    for some $C$.
   We have $X_{2,2}\subset \mathbb P^5$, with a regular $G$-action on the ambient $\mathbb P^5$, so that $X_{2, 2}$ is the $G$-invariant base locus of a pencil of quadrics. Let $\beta : Y \to \mathbb P^5$ be the $G$-equivariant blowup along $X_{2, 2}$, it admits a $G$-equivariant quadric fibration $\pi : Y \to \mathbb P^1$. Let $\mathcal W \to \mathbb P^1$ be the relative Fano variety of planes of $Y/\mathbb P^1$, 
   with Stein factorization
   $$
   \mathcal W \to C\to \mathbb P^1,
   $$
   here $C$ is a smooth projective curve of genus $2$ with a regular $G$-action.
   The proof of \cite[Proposition 7.4]{LT19} shows that $\mathcal W=\cF_4(X_{2,2})$, the space of plane conics on $X_{2, 2}$. We have
    the Abel-Jacobi map
    \[
    \mathrm{AJ} : \mathcal W \to (\mathrm{\mathbf{CH}}^2_{X_{2, 2}/\mathbb C})^{2\ell}.
    \]
    The fibers of $\mathcal W \to C$ are orthogonal Grassmannians, and thus rational. This shows that the Abel-Jacobi map factors as
    \[
    \mathcal W \to C \to (\mathrm{\mathbf{CH}}^2_{X_{2, 2}/\mathbb C})^{2\ell}.
    \]
    This induces a $G$-equivariant isomorphism
    \[
    \Pic^1(C) \cong (\mathrm{\mathbf{CH}}^2_{X_{2, 2}/\mathbb C})^{2\ell}.
    \]
\end{proof}

Now, assume that the $G$-action on $X_{2, 2}$ is projectively linearizable. By Lemma~\ref{lemma:X22}, 
\[
2[\cF_2(X_{2, 2})] = [\mathrm{Pic}^1(C)].
\]
On the other hand, by  Theorem~\ref{theorem:IJTobstructions}, there exists a smooth projective $G$-curve $C'$ such that 
\[
[\cF_2(X_{2, 2})] = [\Pic^e(C')],
\]
for some $e$. By Torelli, $C'$ is $G$-equivariantly isomorphic to $C$. Thus, 
\[
0 =2e[\Pic^1(C)] = [\Pic^1(C)],
\]
because $C$ has genus $2$.
This implies that $\cF_2(X_{2, 2})$ is a trivial $\IJ(X_{2, 2})$-torsor, as claimed.

\subsection*{Prime Fano threefolds of degree 16}
Consider a $G$-Fano threefold $X_{16}\subset \bP^{10}$  
with projectively linearizable $G$-action.    
Then there exists a smooth projective curve $C$, of genus $3$,  with a regular $G$-action such that
\[
[(\mathrm{\mathbf{CH}}^2_{X_{16}/\mathbb C})^\ell ] = [\Pic^e(C)],
\]
for some $e$ where $\ell$ is the class of lines. The variety of conics
$\cF_2(X_{16})$ is a $\bP^1$-bundle over $C$, see, e.g., 
 \cite[Theorem 9.2]{KuzP}.
It follows that 
\[
2e[\Pic^1(C)] = [(\mathrm{\mathbf{CH}}^2_{X_{16}/\mathbb C})^{2\ell}] = [\Pic^1(C)].
\]
Indeed the last equality follows from \cite[Corollary 8.9]{KuzP}.
On the other hand, since $D$ has genus $3$, we have
\[
4[\Pic^1(C)] = 0.
\]
Together, these imply 
that $[(\mathrm{\mathbf{CH}}^2_{X_{16}/\mathbb C})^{3\ell} ] = 0$.
By \cite[Theorem 9.4]{KuzP}, the Abel-Jacobi map induces a $G$-equivariant isomorphism
\[
\mathrm{AJ} : \cF_3(X_{16}) \to (\mathrm{\mathbf{CH}}^2_{X_{16}/\mathbb C})^{3\ell}, 
\]
and the assertion of Theorem~\ref{thm:fanolinear} follows.

Over nonclosed fields, we have the converse by \cite{KuzP}: 
the existence of a twisted cubic curve over the ground field 
implies birationality to (a form of) $V_5$, and these are rational, over any ground field. 

In the equivariant context, we follow the arguments of \cite{KuzP}. Suppose that we have a $G$-stable twisted cubic $R$. Assume $R$ does not admit a bisecant line. Then the Sarkisov link with center at $R$ yields an equivariant birationality to $V_5$ (see \cite[Theorem 5.10]{KuzP}). However, this does not immediately suffice to conclude projective linearizability. Recall that $\Aut(V_5)=\PGL_2$. We have already mentioned that $V_5$ is rigid for the action of $\fA_5\subset \PGL_2$, by \cite{CS}. 
However, $\fA_5$ cannot act on $X_{16}$, since it does not act 
on the Jacobian of a genus 3 curve. 
Linearizations of cyclic and dihedral group actions on $V_5$ are given in Example 4.34, respectively, Example 4.37 of \cite{C-Calabi}. The remaining case of a possible $\fS_4$-action on $V_5$ does not allow an equivariant birationality to $X_{16}$; indeed this would involve blowing up an $\fS_4$-invariant genus 3 and degree 9 curve $T \subset V_5$, which would intersect the $\fS_4$-invariant hyperplane section of $V_5$ in 9 points, which is impossible (we are grateful to I. Cheltsov for providing us with the proof and the reference). When $R$ admits a bisecant line, the intersection of the linear span of $R$ with $X_{16}$ is the union of $R$ and a bisecant line $\ell$ which is $G$-stable. Then equivariant birationality to $\bP^3$ can be realized using the Sarkisov link with center at $\ell$ (see \cite[Theorem 5.8]{KuzP}).

In other cases, according to \cite[Lemma 2.9]{KuzP}, either we have a $G$-stable line or there is a $G$-stable union of three lines. When there is a $G$-stable line, one may argue as above. Finally when we have a $G$-stable union of three lines, they meet at a single point $x$ which is $G$-fixed. By \cite[Lemma 5.11 and Lemma 5.16]{KuzP}, the double projection from $x$ yields equivariant birationality to a complete intersection of three quadrics $\overline{X}_{16}$ containing a $G$-stable plane $\Pi$. 
Projection from $\Pi$ provides an equivariant birational map
\[
\overline{X}_{16}\dashrightarrow \mathbb P^3, 
\]
see, e.g., \cite[Corollary 5.15]{KuzP}.
Thus $X_{16}$  is projectively linearizable.

\subsection*{Prime Fano threefolds of degree 18}
Let $X_{18}\subset \bP^{11}$ be a projectively linearizable $G$-Fano threefold.
Then there is a smooth projective genus $2$ curve $C$ with a regular $G$-action such that
\[
[(\mathrm{\mathbf{CH}}^2_{X_{18}/\mathbb C})^\ell ] = [\Pic^e(C)],
\]
for some $e$, and therefore 
\[
[(\mathrm{\mathbf{CH}}^2_{X_{18}/\mathbb C})^{2\ell} ] =  e[\Pic^2(C)] =0.
\]
By \cite[Corollary 9.10]{KuzP}, the Abel-Jacobi map gives a $G$-equivariant isomorphism
\[
\mathrm{AJ} : \cF_2(X_{18}) \cong (\mathrm{\mathbf{CH}}^2_{X_{18}/\mathbb C})^{2\ell},
\]
and thus $\cF_2(X_{18})^G\neq \emptyset$.

Again, over a nonclosed field $K$, we have a converse statement (\cite{KuzP}), which requires, in addition to a conic $R\subset X_{18}$ over $K$, a (sufficiently general) rational point $x\in X_{18}(K)$. These yield birationality of $X_{18}$ to a sextic del Pezzo fibration $\varphi: \cX\to \bP^1$, together with a bisection, coming from $R$, and a trisection, from $x$, over $K$. 
Sextic del Pezzo surfaces are rational over a nonclosed field, provided they have points over degree 2 and degree 3 extensions; thus,  
$\cX$ is rational over $K(\bP^1)$, and $X_{18}$ is rational over $K$.     
A similar construction in the equivariant context
would require, in addition to a $G$-invariant conic, a $G$-fixed point. However,  
even a section of a sextic del Pezzo fibration  $\varphi$ does not guarantee linearizability, see \cite[Section 5, Case (F)]{CTZ-segre}!

It would be interesting to settle the (projective) linearizability of actions on $X_{18}$.

\bibliographystyle{alpha} 
\bibliography{Bibliography.bib}

\newcommand{\etalchar}[1]{$^{#1}$}
\begin{thebibliography}{ACC{\etalchar{+}}23}

\bibitem[ACC{\etalchar{+}}23]{C-Calabi}
C.~Araujo, A.-M. Castravet, I.~Cheltsov, K.~Fujita, A.-S. Kaloghiros,
  J.~Martinez-Garcia, C.~Shramov, H.~S\"{u}\ss, and N.~Viswanathan.
\newblock {\em The {C}alabi problem for {F}ano threefolds}, volume 485 of {\em
  London Mathematical Society Lecture Note Series}.
\newblock Cambridge University Press, Cambridge, 2023.

\bibitem[AT19]{AT19}
D.~Abramovich and M.~Temkin.
\newblock Functorial factorization of birational maps for qe schemes in
  characteristic 0.
\newblock {\em Algebra Number Theory}, 13(2):379--424, 2019.

\bibitem[Bea77]{Beau77}
A.~Beauville.
\newblock Vari\'{e}t\'{e}s de {P}rym et jacobiennes interm\'{e}diaires.
\newblock {\em Ann. Sci. \'{E}cole Norm. Sup. (4)}, 10(3):309--391, 1977.

\bibitem[Bla07]{blanc2}
J.~Blanc.
\newblock The number of conjugacy classes of elements of the {C}remona group of
  some given finite order.
\newblock {\em Bull. Soc. Math. France}, 135(3):419--434, 2007.

\bibitem[BLRT23]{BLRT23}
R.~Beheshti, B.~Lehmann, E.~Riedl, and S.~Tanimoto.
\newblock Rational curves on del {P}ezzo surfaces in positive characteristic.
\newblock {\em Trans. Amer. Math. Soc. Ser. B}, 10:407--451, 2023.

\bibitem[Bru08]{Bruin}
N.~Bruin.
\newblock The arithmetic of {P}rym varieties in genus 3.
\newblock {\em Compositio Mathematica}, 144(2):317--338, mar 2008.

\bibitem[BS83]{BS}
S.~Bloch and V.~Srinivas.
\newblock Remarks on correspondences and algebraic cycles.
\newblock {\em Amer. J. Math.}, 105(5):1235--1253, 1983.

\bibitem[BvBT24]{BBT}
Chr. B\"ohning, H.-Chr.~Graf von Bothmer, and Yu. Tschinkel.
\newblock Equivariant birational types and derived categories.
\newblock {\em Math. Nachr.}, 297(11):4333--4355, 2024.

\bibitem[BW23]{BW19}
O.~Benoist and O.~Wittenberg.
\newblock Intermediate {J}acobians and rationality over arbitrary fields.
\newblock {\em Ann. Sci. \'Ec. Norm. Sup\'er.}, 56(4):1029--1084, 2023.

\bibitem[CS16]{CS}
I.~Cheltsov and C.~Shramov.
\newblock {\em Cremona groups and the icosahedron}.
\newblock Monogr. Res. Notes Math. Boca Raton, FL: CRC Press, 2016.

\bibitem[CTZ23]{CTZ-segre}
I.~Cheltsov, Yu. Tschinkel, and Zh. Zhang.
\newblock Equivariant geometry of the {S}egre cubic and the {B}urkhardt
  quartic.
\newblock {\tt arXiv:2308.15271}, 2023.

\bibitem[FJS{\etalchar{+}}23]{FJSV22}
S.~Frei, L.~Ji, S.~Sankar, Viray B., and I.~Vogt.
\newblock Curve classes on conic bundle threefolds and applications to
  rationality.
\newblock {\em Algebraic Geometry}, 2023.
\newblock to appear, {\tt arXiv:2207.07093}.

\bibitem[Ful98]{Fulton}
W.~Fulton.
\newblock {\em Intersection theory}.
\newblock Springer-Verlag, Berlin, second edition, 1998.

\bibitem[GHS03]{GHS}
T.~Graber, J.~Harris, and J.~Starr.
\newblock Families of rationally connected varieties.
\newblock {\em J. Amer. Math. Soc.}, 16(1):57--67, 2003.

\bibitem[Har77]{Hartshorne}
R.~Hartshorne.
\newblock {\em Algebraic geometry}, volume No. 52 of {\em Graduate Texts in
  Mathematics}.
\newblock Springer-Verlag, New York-Heidelberg, 1977.

\bibitem[HT12]{HT12}
B.~Hassett and Yu. Tschinkel.
\newblock Spaces of sections of quadric surface fibrations over curves.
\newblock In {\em Compact moduli spaces and vector bundles}, volume 564 of {\em
  Contemp. Math.}, pages 227--249. Amer. Math. Soc., Providence, RI, 2012.

\bibitem[HT21a]{HT-cycle}
B.~Hassett and Yu. Tschinkel.
\newblock Cycle class maps and birational invariants.
\newblock {\em Commun. Pure Appl. Math.}, 74(12):2675--2698, 2021.

\bibitem[HT21b]{HT-quad}
B.~Hassett and Yu. Tschinkel.
\newblock Rationality of complete intersections of two quadrics over nonclosed
  fields.
\newblock {\em Enseign. Math.}, 67(1-2):1--44, 2021.
\newblock With an appendix by Jean-Louis Colliot-Th\'{e}l\`ene.

\bibitem[HT21c]{HT-18}
B.~Hassett and Yu. Tschinkel.
\newblock Rationality of {F}ano threefolds of degree 18 over non-closed fields.
\newblock In {\em Rationality of varieties}, volume 342 of {\em Progr. Math.},
  pages 237--247. Birkh\"{a}user/Springer, Cham, [2021] \copyright 2021.

\bibitem[HT22]{HT-odd}
B.~Hassett and Yu. Tschinkel.
\newblock Equivariant geometry of odd-dimensional complete intersections of two
  quadrics.
\newblock {\em Pure Appl. Math. Q.}, 18(4):1555--1597, 2022.

\bibitem[HT23]{HT23}
B.~Hassett and Yu. Tschinkel.
\newblock Torsors and stable equivariant birational geometry.
\newblock {\em Nagoya Math. J.}, 250:275--297, 2023.

\bibitem[JJ23]{JJ23}
L.~Ji and M.~Ji.
\newblock Rationality of real conic bundles with quartic discriminant curve.
\newblock {\em Int. Math. Res. Not. IMRN}, 2023.
\newblock Online publication.

\bibitem[Kah21]{Kahn}
B.~Kahn.
\newblock On the universal regular homomorphism in codimension 2.
\newblock {\em Ann. Inst. Fourier (Grenoble)}, 71(2):843--848, 2021.

\bibitem[KP23]{KuzP}
A.~Kuznetsov and Yu. Prokhorov.
\newblock Rationality of {Fano} threefolds over non-closed fields.
\newblock {\em Am. J. Math.}, 145(2):335--411, 2023.

\bibitem[KT22]{BnG}
A.~Kresch and Yu. Tschinkel.
\newblock Equivariant birational types and {B}urnside volume.
\newblock {\em Ann. Sc. Norm. Super. Pisa Cl. Sci. (5)}, 23(2):1013--1052,
  2022.

\bibitem[LRT23]{LRT23}
B~Lehmann, E.~Riedl, and S.~Tanimoto.
\newblock Non-free sections of {F}ano fibrations, 2023.
\newblock {\tt arXiv:2301.01695}.

\bibitem[LT19]{LT19}
B.~Lehmann and S.~Tanimoto.
\newblock Geometric {M}anin's conjecture and rational curves.
\newblock {\em Compos. Math.}, 155(5):833--862, 2019.

\bibitem[LT24]{LT24}
B.~Lehmann and S.~Tanimoto.
\newblock Classifying sections of del {P}ezzo fibrations, {I}.
\newblock {\em J. Eur. Math. Soc. (JEMS)}, 26(1):289--354, 2024.

\bibitem[Mat95]{Matsuki}
K.~Matsuki.
\newblock Weyl groups and birational transformations among minimal models.
\newblock {\em Mem. Amer. Math. Soc.}, 116(557):vi+133, 1995.

\bibitem[Mil08]{Milne}
J.~Milne.
\newblock Abelian {V}arieties, 2008.
\newblock online lecture notes.

\bibitem[MM81]{Mori1981}
Sh. Mori and Sh. Mukai.
\newblock Classification of {F}ano $3$-folds with {$B_2 \geq 2$}.
\newblock {\em Manuscripta mathematica}, 36:147--162, 1981.

\bibitem[MM83]{MM83}
Sh. Mori and Sh. Mukai.
\newblock On {F}ano {$3$}-folds with {$B\sb{2}\geq 2$}.
\newblock In {\em Algebraic varieties and analytic varieties ({T}okyo, 1981)},
  volume~1 of {\em Adv. Stud. Pure Math.}, pages 101--129. North-Holland,
  Amsterdam, 1983.

\bibitem[MM03]{MM03}
Sh. Mori and Sh. Mukai.
\newblock Erratum: ``{C}lassification of {F}ano 3-folds with {$B_2\geq 2$}''
  [{M}anuscripta {M}ath. {\bf 36} (1981/82), no. 2, 147--162; {MR}0641971
  (83f:14032)].
\newblock {\em Manuscripta Math.}, 110(3):407, 2003.

\bibitem[Mum74]{prym}
D.~Mumford.
\newblock Prym varieties. {I}.
\newblock In {\em Contributions to analysis (a collection of papers dedicated
  to {L}ipman {B}ers)}, pages 325--350. 1974.

\bibitem[Mur85]{Murre}
J.~P. Murre.
\newblock Applications of algebraic {$K$}-theory to the theory of algebraic
  cycles.
\newblock In {\em Algebraic geometry, {S}itges ({B}arcelona), 1983}, volume
  1124 of {\em Lecture Notes in Math.}, pages 216--261. Springer, Berlin, 1985.

\bibitem[Pro13]{pro-inv}
Yu.~G. Prokhorov.
\newblock On birational involutions of {{\(\mathbb P^3\)}}.
\newblock {\em Izv. Math.}, 77(3):627--648, 2013.

\bibitem[Pro18]{Prokhorov2018}
Yu.~G. Prokhorov.
\newblock The rationality problem for conic bundles.
\newblock {\em Russian Mathematical Surveys}, 73(3):375--456, 2018.

\bibitem[Sca21]{Scavia}
F.~Scavia.
\newblock Equivariant intermediate {J}acobians and intersections of two
  quadrics, 2021.
\newblock unpublished note.

\bibitem[Tes09]{Testa}
D.~Testa.
\newblock The irreducibility of the spaces of rational curves on del {P}ezzo
  surfaces.
\newblock {\em J. Algebraic Geom.}, 18(1):37--61, 2009.

\bibitem[TYZ24]{TYZ23}
Yu. Tschinkel, K.~Yang, and Zh. Zhang.
\newblock Equivariant birational geometry of linear actions.
\newblock {\em EMS Surv. Math. Sci.}, 11(2):235--276, 2024.

\end{thebibliography}

\end{document}